\documentclass[A4j,11pt]{article}
\usepackage{graphics}
\usepackage{amsmath}
\usepackage{amssymb}
\usepackage{latexsym}
\usepackage{amsfonts}
\begin{document}

\title{{\bf Hermann type actions on\\
a pseudo-Riemannian symmetric space}}
\author{{\bf Naoyuki Koike}}
%
%
\date{}
%
\maketitle
\begin{abstract}
In this paper, we first investigate the geometry of the orbits of the isotropy 
action of a semi-simple pseudo-Riemannian symmetric space by investigating 
the complexified action.  Next we investigate the geometry of the orbits of a 
Hermann type action on a semi-simple pseudo-Riemannian symmetric space.  
By considering two special Hermann type actions on a semi-simple 
pseudo-Riemannian symmetric space, we recognize an interesting structure of 
the space.  As a special case, we we recognize an interesting structure of 
the complexification of a semi-simple pseudo-Riemannian symmetric space.  
Also, we investigate a homogeneous submanifold with flat section in a 
pseudo-Riemannian symmmtric space under certain conditions.  
\end{abstract}

\section{Introduction}
In Riemannian symmetric spaces, the notion of an equifocal submanifold was 
introduced by Terng-Thorbergsson in [36].  This notion is defined as a compact submanifold with flat section 
such that the normal holonomy group is trivial and that the focal radius 
functions for each parallel normal vector field are constant.  
However, the condition of the equifocality is rather weak in the case 
where the Riemannian symmetric spaces are of non-compact type and the 
submanifold is non-compact.  So we [17,18] 
have recently introduced the notion of a complex equifocal 
submanifold in a Riemannian symmetric space $G/K$ of non-compact type.  
This notion is defined by imposing the constancy of the complex focal 
radius functions in more general.  
Here we note that the complex focal radii are the 
quantities indicating the positions of the focal points of the 
extrinsic complexification of the submanifold, where the submanifold needs 
to be assumed to be complete and of class $C^{\omega}$ 
(i.e., real analytic).  
On the other hand, Heintze-Liu-Olmos [13] has recently defined the 
notion of an isoparametric submanifold with flat section in a general 
Riemannian manifold as a submanifold such that the normal holonomy group is 
trivial, its sufficiently close parallel submanifolds are of constant mean 
curvature with respect to the radial direction and that the image of the 
normal space at each point by the normal 
exponential map is flat and totally geodesic.  
We [18] showed the following fact:

\vspace{0.2truecm}

{\sl All isoparametric submanifolds with 
flat section in a Riemannian symmetric space $G/K$ of non-compact type are 
complex equifocal and that conversely, all curvature-adapted and complex 
equifocal submanifolds are isoparametric ones with flat section.}

\vspace{0.2truecm}

\noindent
Here the curvature-adaptedness means that, for each normal vector $v$ of 
the submanifold, the Jacobi operator $R(\cdot,v)v$ preserves the tangent space 
of the submanifold invariantly and the restriction of $R(\cdot,v)v$ to the 
tangent space commutes with the shape operator $A_v$, where $R$ is the 
curvature tensor of $G/K$.  
Note that curvature-adapted hypersurfaces in a complex hyperbolic space 
(and a complex projective space) mean so-called Hopf hypersurfaces and that 
curvature-adapted complex equifocal hypersurfaces in the space mean Hopf 
hypersurfaces with constant principal curvatures, which are classified by 
J. Berndt [2].  Also, he [3] classified curvature-adapted hypersurfaces with 
constant principal curvatures (i.e., curvature-adpated complex equifocal 
hypersurfaces) in the quaternionic hyperbolic space.  
As a subclass of the class of complex equifocal submanifolds, 
we [19] defined the notion of a proper complex equifocal submanifold in 
$G/K$ as a complex equifocal submanifold whose lifted submanifold to 
$H^0([0,1],\mathfrak g)$ ($\mathfrak g:={\rm Lie}\,G$) through 
some pseudo-Riemannian submersion of $H^0([0,1],\mathfrak g)$ onto $G/K$ is 
proper complex isoparametric in the sense of [17], where we note that 
$H^0([0,1],\mathfrak g)$ is a pseudo-Hilbert space consisting of certain kind 
of paths in the Lie algebra $\mathfrak g$ of $G$.  
For a $C^{\omega}$-submanifold $M$, we [18] showed that $M$ is proper 
complex equifocal if and only if 
the lift of the complexification $M^{\bf c}$ (which is a submanifold 
in the anti-Kaehlerian symmetric space $G^{\bf c}/K^{\bf c}$) of $M$ 
to $H^0([0,1],\mathfrak g^{\bf c})$ 
($\mathfrak g^{\bf c}:={\rm Lie}\,G^{\bf c}$) 
by some anti-Kaehlerian submersion of $H^0([0,1],\mathfrak g^{\bf c})$ onto 
$G^{\bf c}/K^{\bf c}$ is proper anti-Kaehlerian isoparametric in the sense of 
[18].  
This fact implies that a proper complex equifocal submanifold is a complex 
equifocal submanifold whose complexification has regular focal structure.   
Let $G/K$ be a Riemannian symmetric space of non-compact type 
and $H$ be a closed subgroup of $G$.  If the $H$-action is proper and 
there exists a complete embedded flat submanifold meeting all $H$-orbits 
orthogonally, then it is called a {\it complex hyperpolar action}.  
Principal orbits of a complex hyperpolar action are complex equifocal.  
If $H$ is a symmetric subgroup of $G$ (i.e., 
$({\rm Fix}\,\sigma)_0\subset H\subset{\rm Fix}\,\sigma$ for some involution 
$\sigma$ of $G$), then the $H$-action is called a {\it Hermann type action}, 
where ${\rm Fix}\,\sigma$ is the fixed point group of $\sigma$ and 
$({\rm Fix}\,\sigma)_0$ is the identity component of the group.  
Hermann type actions are complex hyperpolar.  
We ([18,19]) showed the following fact:

\vspace{0.2truecm}

{\sl Principal orbits of a Hermann type action are curvature-adapted 
and proper complex equifocal.}

\vspace{0.2truecm}

\noindent
Similarly, we can define the notions of a complex equifocal submanifold, 
proper complex equifocal one and a curvature-adapted one in a 
pseudo-Riemannian symmetric space (see Section 2).  
Also, we can define the notions of a complex hyperpolar action and a Hermann 
type action on a pseudo-Riemannian symmetric space.  
We [23] showed the following fact:

\vspace{0.2truecm}

{\sl All isoparametric submanifolds with flat section in a pseudo-Riemannian 
symmetric space $G/K$ are complex equifocal.  Conversely all 
curvature-adapted complex equifocal submanifolds such that $A$ and 
$R(\cdot,v)v$ are semi-simple for any normal vector $v$ 
are isoparametric ones with flat section, where $A_v$ is the 
shape operator and $R$ is the curvature tensor of $G/K$ and the 
semi-simplenesses of $A_v$ and $R(\cdot,v)v$ mean that their complexifications 
are diagonalizable.}

\vspace{0.2truecm}

L. Geatti and C. Gorodski [9] has recently showed that a polar representation 
of a real reductive algebraic group on a pseudo-Eucliean space has the same 
closed orbits as the isotropy representation (i.e., the linear isotropy 
action) of a pseudo-Riemannian symmetric space (see Theorem 1 of [9]).  
Also, they showed that the principal orbits of the polar representation 
through a semi-simple element (i.e., the orbit through a regular element 
(in the sense of [9])) is an isoparametric submanifold by investigating 
the complexified representation (see Theorem 11 (also Example 12) of [9]), 
where an isoparametric submanifold means a pseudo-Riemannian submanifold such 
that the (restricted) normal holonomy group is trivial and that 
the shape operator for each (local) parallel normal vector field is 
semi-simple and has constant complex principal curvature.  
All isoparametric submanifold in this sense are isoparametric ones 
(with flat section) in the sense of [13].   
Let $G/H$ be a (semi-simple) pseudo-Riemannian symmetric space 
(equipped with the metric $\langle\,\,,\,\,\rangle$ induced from the Killing 
form of the Lie algebra $\mathfrak g$ of $G$).  
In this paper, we first investigate the complexified shape operators of 
the orbits of the isotropy action of $G/H$ (i.e., the $H$-action on $G/H$) 
by investigating the orbits of the isotropy action of $G^{\bf c}/H^{\bf c}$ 
(see Section 3).  
Next, by using the investigation, we prove the following fact for 
the orbits of Hermann type action.  

\vspace{0.5truecm}

\noindent
{\bf Theorem A.} {\sl Let $G/H$ be a (semi-simple) pseudo-Riemannian 
symmetric space, 
$H'$ be a symmetric 
subgroup of $G$, $\sigma$ (resp. $\sigma'$) be an involution of $G$ 
with $({\rm Fix}\,\sigma)_0\subset H\subset{\rm Fix}\,\sigma$ (resp. 
$({\rm Fix}\,\sigma')_0\subset H'\subset{\rm Fix}\,\sigma'$), 
$L:=({\rm Fix}(\sigma\circ\sigma'))_0$ and $\mathfrak l:={\rm Lie}\,L$.  
Assume that $G$ is not compact and 
$\sigma\circ \sigma'=\sigma'\circ \sigma$.  Then the following statements 
${\rm(i)}$ and ${\rm(ii)}$ hold:

{\rm (i)} The orbit $H'(eH)$ of the $H'$-action on $G/H$ is a reflective 
pseudo-Riemannian submanifold and it is homothetic to the semi-simple 
pseudo-Riemannian symmetric space $H'/H\cap H'$.  
For each $x\in H'(eH)$, the section $\Sigma_x$ of $H'(eH)$ through $x$ is 
homothetic to the semi-simple pseudo-Riemannian symmetric space $L/H\cap H'$.  


{\rm (ii)} Let $M$ be a principal orbit 
of the $H'$-action through a point $\exp_G(w)H$ ($w\in \mathfrak q\cap
\mathfrak q'$ s.t. ${\rm ad}(w)\vert_{\mathfrak l}\,:\,$semi-simple) of 
$\Sigma_{eH}\setminus F$, where 
$\mathfrak q:={\rm Ker}(\sigma+{\rm id})(=T_{eH}(G/H))$, 
$\mathfrak q':={\rm Ker}(\sigma'+{\rm id})$ and 
$F$ is a focal set of $H'(eH)$.  Then $M$ is curvature-adapted 
and proper complex equifocal, for any normal vector $v$ of $M$, 
$R(\cdot,v)v$ and the shape operator $A_v$ are semi-simple.  
Hence it is an isoparametric submanifold with flat section.}

\vspace{0.5truecm}

\noindent
{\it Remark 1.1.} (i) Since 
$\displaystyle{\mathop{\cup}_{w\in\mathfrak q\cap\mathfrak q'\,\,{\rm s.t.}\,\,
{\rm ad}(w)\vert_{\mathfrak l}:{\rm semi-simple}}(H'\cap H)(\exp_G(w)H)}$
is an open dense subset of $L(eH)$, it is shown that 
$$\mathop{\cup}_{w\in\mathfrak q\cap\mathfrak q'\,\,{\rm s.t.}\,\,
{\rm ad}(w)\vert_{\mathfrak l}:{\rm semi-simple}}H'(\exp_G(w)H)$$
is an open dense subset of $G/H$.  

(ii) It is shown that, if $M$ is a curvature-adapted complex equifocal 
submanifold and, for any normal vector $v$ of $M$, $R(\cdot,v)v$ 
and $A_v$ are semi-simple, then it is an isoparametric submanifold 
with flat section (see Proposition 9.1 of [23]).  

(iii) When we take a Riemannian symmetric space of 
non-compact type as $G/H$ in this theorem, we have 
$\displaystyle{\mathop{\cup}_{x\in H'(eH)}\Sigma_x=G/H}$ and $F=\emptyset$.  

\vspace{0.5truecm}

L. Geatti [8] has recently defined a pseudo-Kaehlerian structure on some 
$G$-invariant domain of the complexification $G^{\bf c}/H^{\bf c}$ 
of a semi-simple pseudo-Riemannian symmetric space $G/H$.  
On the other hand, 
we [23] have recently defined an anti-Kaehlerian structure on the whole of 
the complexification $G^{\bf c}/H^{\bf c}$.  
By applying Theorem A to the complexification $G^{\bf c}/H^{\bf c}$ 
(equipped with the natural anti-Kaehlerian structure) of 
a semi-simple pseudo-Riemannian symmetric space $G/H$ and a symmetric subgroup 
$G$ of $G^{\bf c}$, we recognize an interesting structure of 
$G^{\bf c}/H^{\bf c}$.  Here we note that an involution 
$\sigma$ of $G^{\bf c}$ with $({\rm Fix}\,\sigma)_0\subset H^{\bf c}
\subset{\rm Fix}\,\sigma$ and the conjugation $\tau$ of $G^{\bf c}$ with 
respect to $G$ are commutative.  In this case, the group corresponding to $L$ 
in the statement of Theorem A is the dual $G^{\ast_H}$ of $G$ with respect to 
$H$.  Hence we have the following fact.  



\centerline{
\unitlength 0.1in
%
}

\vspace{0.5truecm}

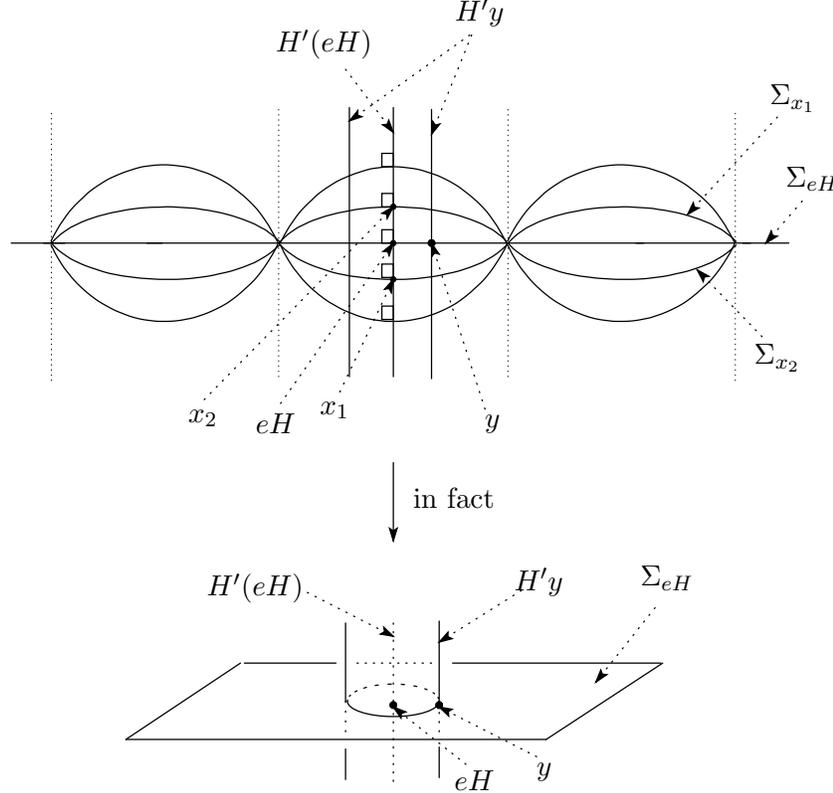
\centerline{{\bf Fig. 1.}}

\vspace{0.5truecm}

\noindent
{\bf Corollary B.} {\sl Let $G^{\bf c}/H^{\bf c}$ and $G^{\ast_H}$ be 
as above.  Then the following statements ${\rm (i)}$ and ${\rm (ii)}$ hold:

{\rm (i)} The orbit $G(eH^{\bf c})$ is a reflective pseudo-Riemannian 
submanifold and it is homothetic to the pseudo-Riemannian symmetric space 
$G/H$.  
For each $x\in G(eH^{\bf c})$, the section $\Sigma_x$ of 
$G(eH^{\bf c})$ through $x$ is homothetic to the pseudo-Riemannian symmetric 
space $G^{\ast_H}/H$.  

{\rm (ii)} For principal orbits of the $G$-action on $G^{\bf c}/H^{\bf c}$, 
the same fact as the statement {\rm (ii)} of Theorem A holds.}

\vspace{0.5truecm}

By considering two special Hermann type actions on a semi-simple 
pseudo-Riemannian symmetric space, we obtain the following interesting fact 
for the structure of the semi-simple pseudo-Riemannian symmetric space.  

\vspace{0.5truecm}

\noindent
{\bf Theorem C.} {\sl Let $G/H$ and $\sigma$ be as in Theorem A, 
$\theta$ the Cartan involution of $G$ with 
$\theta\circ\sigma=\sigma\circ\theta$, $K:=({\rm Fix}\,\theta)_0$ and 
$L:=({\rm Fix}(\sigma\circ\theta))_0$.  
Then the following statements {\rm (i)} and {\rm (ii)} hold:

{\rm (i)} The orbits $K(eH)$ and $L(eH)$ are reflective submanifolds 
satisfying $T_{eH}(G/H)=T_{eH}(K(eH))\oplus T_{eH}(L(eH))$ (orthogonal 
direct sum), $K(eH)$ is anti-homothetic to the Riemannian 
symmetric space $K/H\cap K$ of compact type 
and $L(eH)$ is homothetic to the Riemannian symmetric space $L/H\cap K$ of 
 non-compact type.  
Also, the orbit $K(eH)$ has no focal point.  

{\rm (ii)} Let $M_1$ be a principal orbit of the $K$-action and 
$M_2$ be a principal orbit of the $L$-action through a point of 
$K(eH)\setminus F$, where $F$ is the focal set of $L(eH)$.  
Then $M_i$ ($i=1,2$) are curvature-adapted and proper complex equifocal, 
for any normal vector $v$ of $M_i$, 
$R(\cdot,v)v\vert_{T_xM_i}$ ($x\,:\,$the base point of $v$) and 
the shape operator $A_v$ are diagonalizable.  
Hence they are isoparametric submanifolds with flat section.}

\vspace{0.5truecm}

\noindent
{\it Remark 1.2.} 
For any involution $\sigma$ of $G$, the existence of a Cartan involution 
$\theta$ of $G$ with $\theta\circ\sigma=\sigma\circ\theta$ is assured by 
Lemma 10.2 in [1].  


\vspace{0.5truecm}

\centerline{
\unitlength 0.1in
%
}

\vspace{0.5truecm}

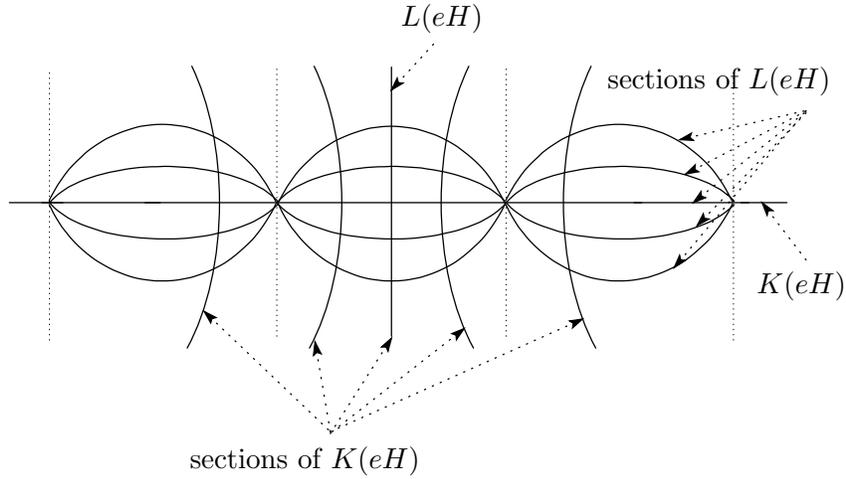
\centerline{{\bf Fig. 2.}}

By applying Theorem C to the complexification $G^{\bf c}/H^{\bf c}$ 
(equipped with the natural anti-Kaehlerian structure) of a 
semi-simple pseudo-Riemannian symmetric space $G/H$, 
we recognize the interesting structure of $G^{\bf c}/H^{\bf c}$.  
In this case, the groups corresponding to $K,\,L$ and $H\cap K$ in the 
statement of Theroem C are as follows.  
Let $\sigma$ be an involution of $G$ with $({\rm Fix}\,\sigma)_0\subset H
\subset{\rm Fix}\,\sigma$, $\theta$ be a Cartan involution of $G$ commuting 
with $\sigma$ and set $K_{\theta}:={\rm Fix}\,\theta$.  Let $G^{\ast}$ be 
the compact dual of $G$ with respect to $K_{\theta}$, 
$H^{\ast}$ be the compact dual of $H$ with respect to $H\cap K_{\theta}$ 
and $(G^d,H^d)$ be the dual of semi-simple symmetric pair $(G,H)$ in the sense 
of [26].  
Then $G^{\ast},\,G^d$ and $H^{\ast}$ correspond to $K,\,L$ and $H\cap K$ 
in the statement of Theorem C, respectively.  
Hence we have the following fact.  

\vspace{0.5truecm}

\noindent
{\bf Corollary D.} {\sl Let $G^{\bf c}/H^{\bf c}$ be 
the complexification (equipped with the natural anti-Kaehlerian structure) of 
a semi-simple pseudo-Riemannian symmetric space $G/H$,  
$G^{\ast}$ (resp. $H^{\ast}$) be the compact dual of $G$ (resp. $H$) and 
$(G^d,H^d)$ be the dual of $(G,H)$.  
Then the following statements {\rm (i)} and {\rm (ii)} hold:

{\rm (i)} The orbits $G^{\ast}(eH^{\bf c})$ and $G^d(eH^{\bf c})$ are 
reflective submanifolds of $G^{\bf c}/H^{\bf c}$ satisfying 
$T_{eH^{\bf c}}(G^{\bf c}/H^{\bf c})=T_{eH^{\bf c}}(G^{\ast}(eH^{\bf c}))
\oplus T_{eH^{\bf c}}(G^d(eH^{\bf c}))$ (orthogonal direct sum), 
$G^{\ast}(eH^{\bf c})$ is anti-homothetic to the Riemannian symmetric space 
$G^{\ast}/H^{\ast}$ of compact type 
and $G^d(eH^{\bf c})$ is homothetic to the Riemannian symmetric space 
$G^d/H^{\ast}$ of non-compact type.  
Also, the orbit $G^{\ast}(eH^{\bf c})$ has no focal point.  

{\rm (ii)} For principal orbits of the $G^{\ast}$-action and $G^d$-action on 
$G^{\bf c}/H^{\bf c}$, the same fact as the statement {\rm(ii)} of Theorem C 
holds.}

\vspace{0.5truecm}

\noindent
{\it Remark 1.3.} In the case where $G/H$ in the statement of Corollary D is 
a Riemannian symmetric space of non-compact type, we have $G^d=G$.  

\vspace{0.5truecm}

Homogeneous submanifolds with flat section in a pseudo-Riemannian symmetric 
space are complex equifocal.  We obtain the following fact for a homogeneous 
submanifold with flat section in a semi-simple pseudo-Riemannian symmetric 
space which admits a reflective focal submanifold, where a reflective 
submanifold means a totally geodesic pseudo-Riemannian submanifold 
with section.  

\vspace{0.5truecm}

\noindent
{\bf Theorem E.} {\sl Let $M$ be a homogeneous submanifold with flat section 
in a semi-simple pseudo-Riemannian symmetric space $G/H$.  Assume that $M$ 
admits a reflective focal submanifold $F$ such that 
$\mathfrak n_{\mathfrak h}(g_{\ast}^{-1}T_{gH}F)$ is a non-degenerate subspace 
of $\mathfrak h$, where $gH$ is an arbitrary point of $F$ and 
$\mathfrak n_{\mathfrak h}(g_{\ast}^{-1}T_{gH}F)$ is the normalizer of 
$g_{\ast}^{-1}T_{gH}F$ in $\mathfrak h$.  Then 
$M$ is a principal orbit of a Hermann type action.  
}

\vspace{0.5truecm}

\noindent
{\it Remark 1.4.} (i) For the $H'$-action in Theorem A, we have 
$\mathfrak n_{\mathfrak h}(T_{eH}(H'(eH)))
=\mathfrak n_{\mathfrak h}(\mathfrak q\cap\mathfrak h')
=\mathfrak h\cap\mathfrak h'+\mathfrak z_{\mathfrak h\cap\mathfrak q'}
(\mathfrak q\cap\mathfrak h')$, where 
$\mathfrak z_{\mathfrak h\cap\mathfrak q'}(\mathfrak q\cap\mathfrak h')$ is 
the centralizer of $\mathfrak q\cap\mathfrak h'$ in 
$\mathfrak h\cap\mathfrak q'$.  
Hence, if $\mathfrak z_{\mathfrak h\cap\mathfrak q'}
(\mathfrak q\cap\mathfrak h')=\{0\}$, then 
$\mathfrak n_{\mathfrak h}(T_{eH}(H'(eH)))$ is a non-degenerate subspace of 
$\mathfrak h$.  Thus almost all principal orbits of the $H'$-action 
have $H'(eH)$ as a reflective focal submanifold as in the statement of 
Theorem E.  

(ii) For the $K$-action in Theorem C, we have 
$\mathfrak n_{\mathfrak h}(T_{eH}(K(eH)))
=\mathfrak n_{\mathfrak h}(\mathfrak q\cap\mathfrak f)
=\mathfrak h\cap\mathfrak f+\mathfrak z_{\mathfrak h\cap\mathfrak p}
(\mathfrak q\cap\mathfrak f)$.  
Hence, $\mathfrak n_{\mathfrak h}(T_{eH}(K(eH)))$ is a non-degenerate 
subspace of $\mathfrak h$.  Similarly, for the $L$-action in Theorem C, 
it is shown that $\mathfrak n_{\mathfrak h}(T_{eH}(L(eH)))$ is a 
non-degenerate subspace of $\mathfrak h$.  
Thus almost all principal orbits of the $K$-action (resp. the $L$-action) 
have $K(eH)$ (resp. $L(eH)$) as a reflective focal submanifold as 
in the statement of Theorem E.  

\section{New notions in a pseudo-Riemannian symmetric space} 
In this section, we shall define new notions in a (semi-simple) 
pseudo-Riemannian symmetric space, which are analogies of notions in 
a Riemannian symmetric space of non-compact type defined in [18].  
Let $M$ be an immersed pseudo-Riemannian submanifold with flat section 
(that, is, $g_{\ast}^{-1}T^{\perp}_xM$ is abelian for any $x=gH\in M$) 
in a (semi-simple) pseudo-Riemannian symmetric space $N=G/H$ (equipped with 
the metric induced from the Killing form of $\mathfrak g:={\rm Lie}\,G$), 
where $T^{\perp}_xM$ is the normal space of $M$ at $x$.  
Denote by $A$ the shape tensor of $M$.  Let $v\in T^{\perp}_xM$ and 
$X\in T_xM$ ($x=gK$), where $T_xM$ is the tangent space of $M$ at $x$.  
Denote by $\gamma_v$ the geodesic in $N$ with $\dot{\gamma}_v(0)=v$, 
where $\dot{\gamma}_v(0)$ is the velocity vector of $\gamma_v$ at $0$.  
The strongly $M$-Jacobi field $Y$ along 
$\gamma_v$ with $Y(0)=X$ (hence $Y'(0)=-A_vX$) is given by 
$$Y(s)=(P_{\gamma_v\vert_{[0,s]}}\circ(D^{co}_{sv}-sD^{si}_{sv}\circ A_v))
(X),\leqno{(2.1)}$$
where $Y'(0)=\widetilde{\nabla}_vY$ ($\widetilde{\nabla}\,:\,$ 
the Levi-Civita connection of $N$), $P_{\gamma_v\vert_{[0,s]}}$ is 
the parallel translation along $\gamma_v\vert_{[0,s]}$ and 
$D^{co}_{sv}$ (resp. $D^{si}_{sv}$) is given by 
$$\begin{array}{c}
\displaystyle{
D^{co}_{sv}=g_{\ast}\circ\cos(\sqrt{-1}{\rm ad}(sg_{\ast}^{-1}v))
\circ g_{\ast}^{-1}}\\
\displaystyle{\left({\rm resp.}\,\,\,\,
D^{si}_{sv}=g_{\ast}\circ
\frac{\sin(\sqrt{-1}{\rm ad}(sg_{\ast}^{-1}v))}
{\sqrt{-1}{\rm ad}(sg_{\ast}^{-1}v)}\circ g_{\ast}^{-1}\right).}
\end{array}$$ 
Here ${\rm ad}$ is the adjoint representation of the Lie algebra 
$\mathfrak g$.  
All focal radii of $M$ along $\gamma_v$ are obtained as real numbers $s_0$ 
with ${\rm Ker}(D^{co}_{s_0v}-s_0D^{si}_{s_0v}\circ A_v)\not=\{0\}$.  So, we 
call a complex number $z_0$ with ${\rm Ker}(D^{co}_{z_0v}-
z_0D^{si}_{z_0v}\circ A_v^{{\bf c}})\not=\{0\}$ a {\it complex 
focal radius of} $M$ {\it along} $\gamma_v$ and call ${\rm dim}\,
{\rm Ker}(D^{co}_{z_0v}-z_0D^{si}_{z_0v}\circ A_v^{{\bf c}})$ the 
{\it multiplicity} of the complex focal radius $z_0$, 
where $A_v^{\bf c}$ is the complexification of $A_v$ and $D^{co}_{z_0v}$ 
(resp. $D^{si}_{z_0v}$) is a ${\bf C}$-linear transformation of 
$(T_xN)^{\bf c}$ defined by 
$$\begin{array}{c}
\displaystyle{
D^{co}_{z_0v}=g^{\bf c}_{\ast}\circ\cos(\sqrt{-1}{\rm ad}^{\bf c}
(z_0g_{\ast}^{-1}v))\circ (g^{\bf c}_{\ast})^{-1}}\\
\displaystyle{\left({\rm resp.}\,\,\,\,
D^{si}_{sv}=g^{\bf c}_{\ast}\circ
\frac{\sin(\sqrt{-1}{\rm ad}^{\bf c}(z_0g_{\ast}^{-1}v))}
{\sqrt{-1}{\rm ad}^{\bf c}(z_0g_{\ast}^{-1}v)}\circ(g^{\bf c}_{\ast})^{-1}
\right),}
\end{array}$$
where $g_{\ast}^{\bf c}$ (resp. ${\rm ad}^{\bf c}$) is the complexification 
of $g_{\ast}$ (resp. ${\rm ad}$).  Here we note that, 
in the case where $M$ is of class $C^{\omega}$, 
complex focal radii along $\gamma_v$ 
indicate the positions of focal points of the (extrinsic) 
complexification $M^{\bf c}(\hookrightarrow G^{\bf c}/H^{\bf c})$ of $M$ 
along the complexified geodesic $\gamma_{\iota_{\ast}v}^{\bf c}$.  Here 
$G^{\bf c}/H^{\bf c}$ is the pseudo-Riemannian symmetric space 
equipped with the metric induced from the Killing form of 
$\mathfrak g^{\bf c}$ regarded as a real Lie algebra 
(which is called the {\it anti-Kaehlerian symmetric space associated with} 
$G/H$), $M^{\bf c}$ and the complexified immersion of $M^{\bf c}$ into 
$G^{\bf c}/H^{\bf c}$ are defined as in [23] and $\iota$ is the natural 
embedding of $G/H$ into $G^{\bf c}/H^{\bf c}$.  
Furthermore, assume that the normal holonomy group of $M$ is trivial.  
Let $\tilde v$ be a parallel unit normal vector field of $M$.  
Assume that the number (which may be $0$ and $\infty$) of distinct complex 
focal radii along $\gamma_{\tilde v_x}$ is independent of the choice of 
$x\in M$.  Furthermore assume that the number is not equal to $0$.  
Let $\{r_{i,x}\,\vert\,i=1,2,\cdots\}$ 
be the set of all complex focal radii along $\gamma_{\tilde v_x}$, where 
$\vert r_{i,x}\vert\,<\,\vert r_{i+1,x}\vert$ or 
"$\vert r_{i,x}\vert=\vert r_{i+1,x}\vert\,\,\&\,\,{\rm Re}\,r_{i,x}
>{\rm Re}\,r_{i+1,x}$" or 
"$\vert r_{i,x}\vert=\vert r_{i+1,x}\vert\,\,\&\,\,
{\rm Re}\,r_{i,x}={\rm Re}\,r_{i+1,x}\,\,\&\,\,
{\rm Im}\,r_{i,x}=-{\rm Im}\,r_{i+1,x}<0$".  
Let $r_i$ ($i=1,2,\cdots$) be complex 
valued functions on $M$ defined by assigning $r_{i,x}$ to each $x\in M$.  
We call these functions $r_i$ ($i=1,2,\cdots$) {\it complex 
focal radius functions for} $\tilde v$.  
We call $r_i\tilde v$ a {\it complex focal normal vector field for} 
$\tilde v$.  If, for each parallel 
unit normal vector field $\tilde v$ of $M$, the number of distinct complex 
focal radii along $\gamma_{\tilde v_x}$ is independent of the choice of 
$x\in M$, each complex focal radius function for $\tilde v$ 
is constant on $M$ and it has constant multiplicity, then 
we call $M$ a {\it complex equifocal submanifold}.  
Also, if parallel submanifolds sufficiently close to $M$ has constant mean 
curvature with respect to the radial direction, then we call $M$ an 
{\it isoparametric submanifold with flat section}.  It is shown that all 
isoparametric submanifolds with flat section are complex equifocal and that, 
conversely, all curvature-adapted complex equifocal submanifold with complex 
diagonalizable shape operators and Jacobi operators are isoparametric 
submanifolds with flat section (see Theorem 9.1 of [23]).  

Let $N=G/H$ be a (semi-simple) pseudo-Riemannian symmetric space and 
$\pi$ be the natural projection of $G$ onto $G/H$.  
Let $\sigma$ be an involution of $G$ with $({\rm Fix}\,\sigma)_0\subset H
\subset{\rm Fix}\,\sigma$ and denote by the same symbol $\sigma$ the 
involution of $\mathfrak g:={\rm Lie}\,G$.  Let $\mathfrak h:=
\{X\in \mathfrak g\,\vert\,\sigma(X)=X\}$ and 
$\mathfrak q:=\{X\in \mathfrak g\,\vert\,\sigma(X)=-X\}$, which is identified 
with the tangent space $T_{eH}N$.  Let $\langle\,\,,\,\,\rangle$ be 
the Killing form of $G$.  Denote by the same symbol 
$\langle\,\,,\,\,\rangle$ both the bi-invariant pseudo-Riemannian metric of 
$G$ induced from $\langle\,\,,\,\,\rangle$ and the pseudo-Riemannian metric of 
$N$ induced from $\langle\,\,,\,\,\rangle$.  Let $\theta$ be a Cartan 
involution of $G$ with $\theta\circ\sigma=\sigma\circ\theta$.  
Denote by the same symbol $\theta$ the involution of $\mathfrak g$ 
induced from $\theta$.  
Let $\mathfrak f:=\{X\in \mathfrak g\,\vert\,\theta(X)=X\}$ and 
$\mathfrak p:=\{X\in \mathfrak g\,\vert\,\theta(X)=-X\}$.  
From $\theta\circ\sigma=\sigma\circ\theta$, it follows that 
$\mathfrak h=\mathfrak h\cap\mathfrak f+\mathfrak h\cap\mathfrak p$ and 
$\mathfrak q=\mathfrak q\cap\mathfrak f+\mathfrak q\cap\mathfrak p$.  
Set $\mathfrak g_+:=\mathfrak p,\,\,\mathfrak g_-:=\mathfrak f$ and 
$\langle\,\,,\,\,\rangle_{\mathfrak g_{\pm}}
:=-\pi_{\mathfrak g_-}^{\ast}\langle\,\,,\,\,\rangle
+\pi_{\mathfrak g_+}^{\ast}\langle\,\,,\,\,\rangle$, where 
$\pi_{\mathfrak g_-}$ (resp. $\pi_{\mathfrak g_+}$) is the projection of 
$\mathfrak g$ onto $\mathfrak g_-$ (resp. $\mathfrak g_+$).  Let 
$H^0([0,1],\mathfrak g)$ be the space of all $L^2$-integrable paths 
$u:[0,1]\to\frak g$ (with respect to 
$\langle\,\,,\,\,\rangle_{\mathfrak g_{\pm}}$).  
It is shown that $(H^0([0,1],\mathfrak g),\,\langle\,\,,\,\,\rangle_0)$ 
is a pseudo-Hilbert space.  Let $H^1([0,1],G)$ be the Hilbert Lie group of all 
absolutely continuous paths $g:[0,1]\to G$ such that the weak derivative $g'$ 
of $g$ is squared integrable (with respect to 
$\langle\,\,,\,\,\rangle_{\mathfrak g_{\pm}}$), that is, 
$g_{\ast}^{-1}g'\in H^0([0,1],\mathfrak g)$.  Define a map 
$\phi:H^0([0,1],\mathfrak g)\to G$ by $\phi(u)=g_u(1)$ 
($u\in H^0([0,1],\mathfrak g)$), where $g_u$ is the element of 
$H^1([0,1],G)$ satisfying $g_u(0)=e$ and $g_{u\ast}^{-1}g_u'=u$.  
We call this map the {\it parallel transport map} (from $0$ to $1$).  
This submersion $\phi$ is a pseudo-Riemannian submersion of 
$(H^0([0,1],\mathfrak g),\langle\,\,,\,\,\rangle_0)$ onto 
$(G,\langle\,\,,\,\,\rangle)$.  
Denote by $\mathfrak g^{\bf c},\mathfrak h^{\bf c},\mathfrak q^{\bf c},
\mathfrak f^{\bf c},\mathfrak p^{\bf c}$ and $\langle\,\,,\,\,\rangle^{\bf c}$ 
the complexifications of $\mathfrak g,\mathfrak h,\mathfrak q,
\mathfrak f,\mathfrak p$ and $\langle\,\,,\,\,\rangle$.  
Set $\mathfrak g^{\bf c}_+:=\sqrt{-1}\mathfrak f+\mathfrak p$ and 
$\mathfrak g^{\bf c}_-:=\mathfrak f+\sqrt{-1}\mathfrak p$.  
Set $\langle\,\,,\,\,\rangle':=2{\rm Re}\langle\,\,,\,\,\rangle^{\bf c}$ and 
$\langle\,\,,\,\,\rangle'_{\mathfrak g^{\bf c}_{\pm}}:=
-\pi^{\ast}_{\mathfrak g^{\bf c}_-}\langle\,\,,\,\,\rangle'+
\pi^{\ast}_{\mathfrak g^{\bf c}_+}\langle\,\,,\,\,\rangle'$, where 
$\pi_{\mathfrak g^{\bf c}_-}$ (resp. $\pi_{\mathfrak g^{\bf c}_+}$) is 
the projection of $\mathfrak g^{\bf c}$ onto $\mathfrak g^{\bf c}_-$ (resp. 
$\mathfrak g^{\bf c}_+$).  
Let $H^0([0,1],\mathfrak g^{\bf c})$ be the space of all $L^2$-integrable 
paths $u:[0,1]\to\mathfrak g^{\bf c}$ (with respect to 
$\langle\,\,,\,\,\rangle'_{\mathfrak g^{\bf c}_{\pm}}$).  
Define a non-degenerate symmetric bilinear form $\langle\,\,,\,\,\rangle'_0$ 
of $H^0([0,1],\mathfrak g^{\bf c})$ by 
$\langle u,v\rangle'_0:=\int_0^1\langle u(t),v(t)\rangle'dt$.  It is shown 
that $(H^0([0,1],\mathfrak g^{\bf c}),\langle\,\,,\,\,\rangle'_0)$ 
is an infinite dimensional anti-Kaehlerian space.  See [18] about the 
definition of an infinite dimensional anti-Kaehlerian space.  
In similar to $\phi$, the parallel transport map $\phi^{\bf c}:H^0([0,1],
\mathfrak g^{\bf c})\to G^{\bf c}$ for $G^{\bf c}$ is defined.  
This submersion $\phi^{\bf c}$ is an anti-Kaehlerian submersion.  
Let $\pi:G\to G/H$ and $\pi^{\bf c}:G^{\bf c}\to G^{\bf c}/H^{\bf c}$ 
be the natural projections.  
By imitating the proof of Theorem 1 of [18], we can show that, 
in the case where $M$ is of class $C^{\omega}$, the following statements 
${\rm(i)}\sim{\rm(iii)}$ are equivalent:

\vspace{0.2truecm}

(i) $M$ is complex equifocal,

(ii) each component of $(\pi\circ\phi)^{-1}(M)$ is complex isoparametric,

(iii) each component of $(\pi^{\bf c}\circ\phi^{\bf c})^{-1}(M^{\bf c})$ 
is anti-Kaehlerian isoparametric.

\vspace{0.2truecm}

\noindent
See [18] about the definitions of a complex isoparametric submanifold and 
an anti-Kaehlerian isoparametric submanifold.  
In particular, if each component of $(\pi\circ\phi)^{-1}(M)$ 
is proper complex isoparametric in the sense of [17], that is, for each normal 
vector $v$ of $(\pi\circ\phi)^{-1}(M)$, there exists a pseudo-orthonormal base 
of the complexified tangent space consisting of the eigenvectors of the 
complexified shape operator for $v$, then we call $M$ a 
{\it proper complex equifocal submanifold}.  
For $C^{\omega}$-submanifold $M$ in $G/H$, it is shown that $M$ is proper 
complex equifocal if and only if $(\pi^{\bf c}\circ\phi^{\bf c})^{-1}
(M^{\bf c})$ is proper anti-Kaehlerian isoparametric in the sense of [18], 
that is, for each normal vector $v$ of $(\pi^{\bf c}\circ\phi^{\bf c})^{-1}
(M^{\bf c})$, there exists a $J$-orthonormal base of the tangent space 
consisting $J$-eigenvectors of the shape operator for $v$, where $J$ is the 
complex structure of $(\pi^{\bf c}\circ\phi^{\bf c})^{-1}(M^{\bf c})$.  
See [18] the definitions of $J$-orthonormal base and $J$-eigenvector.  
Proper anti-Kaehlerian isoparametric submanifolds are interpreted as 
ones having regular focal structure among anti-Kaehlerian isoparametric 
submanifolds.  
From this fact, proper complex equifocal submanifolds are interpreted as ones 
whose complexification has regular focal structure among complex equifocal 
submanifolds.  \newline
$\quad$ Next we shall recall the notions of a complex Jacobi field and 
the parallel translation along a holomorphic curve, which are introduced 
in [23], and we state some facts related to these notions.  
These notions and facts will be used in the next section.  
Let $(M,J,g)$ be an anti-Kaehlerian manifold, $\nabla$ (resp. $R$) be the 
Levi-Civita connection (resp. the curvature tensor) of $g$ 
and $\nabla^{\bf c}$ (resp. $R^{\bf c}$) be the complexification of 
$\nabla$ (resp. $R$).  Let $(TM)^{(1,0)}$ be the holomorphic vector bundle 
consisting of complex vectors of $M$ of type $(1,0)$.  Note that the 
restriction of $\nabla^{\bf c}$ to $TM^{(1,0)}$ is a holomorphic connection 
of $TM^{(1,0)}$ (see Theorem 2.2 of [6]).  
For simplicity, assume that $(M,J,g)$ is complete even if 
the discussion of this section is valid without the assumption of 
the completeness of $(M,J,g)$.  
Let $\gamma:{\bf C}\to M$ be 
a complex geodesic, that is, $\gamma(z)=\exp_{\gamma(0)}(({\rm Re}\,z)
\gamma_{\ast}((\frac{\partial}{\partial s})_0)+({\rm Im}\,z)
J_{\gamma(0)}\gamma_{\ast}((\frac{\partial}{\partial s})_0))$, where $(z)$ is 
the complex coordinate of ${\bf C}$ and $s:={\rm Re}\,z$.  
Let $Y:{\bf C}\to (TM)^{(1,0)}$ be a holomorphic vector 
field along $\gamma$.  
That is, $Y$ assigns $Y_z\in(T_{\gamma(z)}M)^{(1,0)}$ to each 
$z\in {\bf C}$ and, for each holomorphic local coordinate 
$(U,(z_1,\cdots,z_n))$ of $M$ with $U\cap\gamma({\bf C})\not=\emptyset,\,
Y_i:\gamma^{-1}(U)\to{\bf C}$ ($i=1,\cdots,n$) defined by 
$Y_z=\sum\limits_{i=1}^nY_i(z)
(\frac{\partial}{\partial z_i})_{\gamma(z)}$ are holomorphic.  If $Y$ 
satisfies $\nabla^{\bf c}_{\gamma_{\ast}(\frac{d}{dz})}
\nabla^{\bf c}_{\gamma_{\ast}(\frac{d}{dz})}Y+R^{\bf c}(Y,
\gamma_{\ast}(\frac{d}{dz}))\gamma_{\ast}(\frac{d}{dz})=0$, then we call $Y$ 
a {\it complex Jacobi field along} $\gamma$.  
Let $\delta:{\bf C}\times D(\varepsilon)\to M$ be a holomorphic two-parameter 
map, where $D(\varepsilon)$ is the $\varepsilon$-disk centered at $0$ in 
${\bf C}$.  Denote by $z$ (resp. $u$) 
the first (resp. second) parameter of $\delta$.  If 
$\delta(\cdot,u_0):{\bf C}\to M$ is a complex geodesic for each $u_0\in 
D(\varepsilon)$, then we call $\delta$ a {\it complex geodesic variation}.  
It is shown that, for a complex geodesic variation $\delta$, 
the complex variational vector field 
$Y:=\delta_{\ast}(\frac{\partial}{\partial u}\vert_{u=0})$ is a complex Jacobi 
field along $\gamma:=\delta(\cdot,0)$.  
A vector field $X$ on $M$ is said to be {\it real holomorphic} if 
the Lie derivation $L_XJ$ of $J$ with respect to $X$ vanishes.  
It is known that $X$ is a real holomorphic vector field if and only if 
the complex vector field $X-\sqrt{-1}JX$ is holomorphic.  
Let $\gamma:{\bf C}\to M$ be a complex geodesic and 
$Y$ be a holomorphic vector field along $\gamma$.  Denote by $Y_{\bf R}$ 
the real part of $Y$.  Then it is shown that $Y$ is a complex Jacobi field 
along $\gamma$ if and only if, for any $z_0\in{\bf C}$, 
$s\mapsto(Y_{\bf R})_{sz_0}$ is a Jacobi field along the geodesic 
$\gamma_{z_0}(\displaystyle{\mathop{\Longleftrightarrow}_{\rm def}}\,\,
\gamma_{z_0}(s):=\gamma(sz_0))$.  
Next we shall recall the notion of the parallel translation along a 
holomorphic curve.  
Let $\alpha:D\to(M,J,g)$ be a holomorphic curve, where 
$D$ is an open set of ${\bf C}$.  Let $Y$ be a holomorphic vector field along 
$\alpha$.  If $\nabla^{\bf c}_{\alpha_{\ast}(\frac{d}{dz})}Y=0$, then we say 
that $Y$ is {\it parallel}.  
Let $\alpha:D\to (M,J,g)$ be a holomorphic curve.  
For $z_0\in D$ and $v\in(T_{\alpha(z_0)}M)^{(1,0)}$, 
there uniquely exists a parallel holomorphic vector 
field $Y$ along $\alpha$ with $Y_{z_0}=v$.  
We denote $Y_{z_1}$ by 
$(P_{\alpha})_{z_0,z_1}(v)$.  It is clear that $(P_{\alpha})_{z_0,z_1}$ is 
a {\bf C}-linear isomorphism of $(T_{\alpha(z_0)}M)^{(1,0)}$ onto 
$(T_{\alpha(z_1)}M)^{(1,0)}$.  We call $(P_{\alpha})_{z_0,z_1}$ 
the {\it parallel translation along} $\alpha$ 
{\it from} $z_0$ {\it to} $z_1$.  
We consider the case where $(M,J,g)$ is an anti-Kaehlerian symmetric space 
$G^{\bf c}/H^{\bf c}$.  
For $v\in(T_{g_0H^{\bf c}}(G^{\bf c}/H^{\bf c})^{\bf c}$, we define 
${\bf C}$-linear transformations $\widehat D^{co}_v$ and $\widehat D^{si}_v$ 
of $(T_{g_0H^{\bf c}}(G^{\bf c}/H^{\bf c}))^{\bf c}$ by 
$\widehat D^{co}_v:=g_{0\ast}^{\bf c}\circ
\cos(\sqrt{-1}{\rm ad}_{\mathfrak g^{\bf c}}^{\bf c}
((g_{0\ast}^{\bf c})^{-1}v))\circ(g_{0\ast}^{\bf c})^{-1}$ and 
$\widehat D^{si}_v:=g_{0\ast}^{\bf c}\circ
\frac{\sin(\sqrt{-1}{\rm ad}_{\mathfrak g^{\bf c}}^{\bf c}
((g_{0\ast}^{\bf c})^{-1}v))}
{\sqrt{-1}{\rm ad}_{\mathfrak g^{\bf c}}^{\bf c}((g_{0\ast}^{\bf c})^{-1}v)}
\circ(g_{0\ast}^{\bf c})^{-1}$, respectively, where 
${\rm ad}_{\mathfrak g^{\bf c}}^{\bf c}$ is the complexification of the 
adjoint representation ${\rm ad}_{\mathfrak g^{\bf c}}$ of 
$\mathfrak g^{\bf c}$.  
Let $Y$ be a holomorphic vector field along $\gamma^{\bf c}_v$.  Define 
$\widehat Y:D\to(T_{g_0K^{\bf c}}(G^{\bf c}/K^{\bf c}))^{(1,0)}$ by 
$\widehat Y_z:=(P_{\gamma^{\bf c}_v})_{z,0}(Y_z)$ ($z\in D$), where 
$D$ is the domain of $\gamma^{\bf c}_v$.  
Then we have 
$$Y_z=(P_{\gamma^{\bf c}_v})_{0,z}\left(\widehat D^{co}_{zv_{(1,0)}}(Y_0)
+z\widehat D^{si}_{zv_{(1,0)}}(\frac{d\widehat Y}{dz}\vert_{z=0})\right).
\leqno{(2.2)}$$

\section{The isotropy action of a pseudo-Riemannian symmetric space}
In this section, we investigate the complexified shape operators of the orbits 
of the isotropy action of a semi-simple pseudo-Riemannian symmetric space 
by investigating the complexified action.  
Let $G/H$ be a (semi-simple) pseudo-Riemannian symmetric space (equipped with 
the metric $\langle\,\,,\,\,\rangle$ induced from the Killing form $B$ of 
$\mathfrak g$) and $\sigma$ be an involution of $G$ with 
$({\rm Fix}\,\sigma)_0\subset H\subset{\rm Fix}\,\sigma$.  Denote by the same 
symbol $\sigma$ the differential of $\sigma$ at $e$.  Let 
$\mathfrak h:={\rm Lie}\,H$ and $\mathfrak q:={\rm Ker}(\sigma+{\rm id})$, 
which is identified with $T_{eH}(G/H)$.  Let $\theta$ be a Cartan involution 
of $G$ with $\theta\circ\sigma=\sigma\circ\theta,\,\mathfrak f:={\rm Ker}
(\theta-{\rm id})$ and $\mathfrak p:={\rm Ker}(\theta+{\rm id})$.  
Let ${\mathfrak g}^{\bf c},\,{\mathfrak h}^{\bf c},\,{\mathfrak q}^{\bf c},
\,{\mathfrak f}^{\bf c},\,{\mathfrak p}^{\bf c}$ and 
$\langle\,\,,\,\,\rangle^{\bf c}$ be the complexifications of 
$\mathfrak g,\,\mathfrak h,\,\mathfrak q,\,\mathfrak f,\,\mathfrak p$ and 
$\langle\,\,,\,\,\rangle$, respectively.  
The complexification $\mathfrak q^{\bf c}$ is identified with 
$T_{eH^{\bf c}}(G^{\bf c}/H^{\bf c})$.  Under this identification, 
$\sqrt{-1}X\in\mathfrak q^{\bf c}$ corresponds to 
$J_{eH^{\bf c}}X\in T_{eH^{\bf c}}(G^{\bf c}/H^{\bf c})$, where $J$ is 
the complex structure of $G^{\bf c}/H^{\bf c}$.  
Give $G^{\bf c}/H^{\bf c}$ the metric (which also is denoted by 
$\langle\,\,,\,\,\rangle$) induced from the Killing form $B_A$ of 
$\mathfrak g^{\bf c}$ regarded as a real Lie algebra.  Note that $B_A$ 
coincides with $2{\rm Re}\,B^{\bf c}$ and $(J,\langle\,\,,\,\,\rangle)$ is 
an anti-Kaehlerian structure of $G^{\bf c}/H^{\bf c}$, where $B^{\bf c}$ is 
the complexification of $B$.  
Let $\mathfrak a$ be a Cartan subspace of $\mathfrak q$ (that is, 
$\mathfrak a$ is a maximal abelian subspace of $\mathfrak q$ 
and each element of $\mathfrak a$ is semi-simple).  
The dimension of $\mathfrak a$ is called the {\it rank} of $G/H$.  
Without loss of generality, we may assume that 
$\mathfrak a=\mathfrak a\cap\mathfrak f+\mathfrak a\cap\mathfrak p$.  
Let ${\mathfrak q}_{\alpha}^{\bf c}:=\{X\in{\mathfrak q}^{\bf c}\,\vert\,
{\rm ad}(a)^2X=\alpha(a)^2X\,\,{\rm for}\,\,{\rm all}\,\,
a\in{\mathfrak a}^{\bf c}\}$ 
and $\mathfrak h^{\bf c}_{\alpha}:=\{X\in\mathfrak h^{\bf c}\,\vert\,
{\rm ad}(a)^2X=\alpha(a)^2X\,\,{\rm for}\,\,{\rm all}\,\,a\in
\mathfrak a^{\bf c}\}$ for each $\alpha\in({\mathfrak a}^{\bf c})^{\ast}$ 
($({\mathfrak a}^{\bf c})^{\ast}\,:\,$
the (${\bf C}$-)dual space of ${\mathfrak a}^{\bf c}$) and 
$\triangle:=\{\alpha\in({\mathfrak a}^{\bf c})^{\ast}\,\vert\,
{\mathfrak q}^{\bf c}_{\alpha}\not=\{0\}\}$.  Then we have 
$${\mathfrak q}^{\bf c}={\mathfrak a}^{\bf c}+\sum_{\alpha\in\triangle_+}
{\mathfrak q}_{\alpha}^{\bf c}\,\,\,\,\,{\rm and}\,\,\,\,\,
{\mathfrak h}^{\bf c}=\mathfrak z_{\mathfrak h^{\bf c}}({\mathfrak a}^{\bf c})
+\sum_{\alpha\in\triangle_+}{\mathfrak h}_{\alpha}^{\bf c}, \leqno{(3.1)}$$
where $\triangle_+(\subset\triangle)$ is the positive root system under some 
lexicographical ordering and $\mathfrak z_{\mathfrak h^{\bf c}}
(\mathfrak a^{\bf c})$ is the centralizer of $\mathfrak a^{\bf c}$ in 
$\mathfrak h^{\bf c}$.  
Let $\widetilde{\mathfrak a}$ be a Cartan subalgebra of $\mathfrak g$ 
containing $\mathfrak a$ and $\mathfrak g^{\bf c}_{\widetilde{\alpha}}:=
\{X\in\mathfrak g^{\bf c}\,\vert\,{\rm ad}(a)X=\widetilde{\alpha}(a)X\,\,
{\rm for}\,\,{\rm all}\,\,a\in\widetilde{\mathfrak a}^{\bf c}\}$ for each 
$\widetilde{\alpha}\in({\widetilde{\mathfrak a}}^{\bf c})^{\ast}$ and 
$\widetilde{\triangle}:=\{\widetilde{\alpha}\in(\widetilde{\mathfrak a}^{\bf c})^{\ast}\,\vert\,\mathfrak g^{\bf c}_{\widetilde{\alpha}}\not=\{0\}\}$.  
Then we have $\mathfrak g^{\bf c}=\widetilde{\mathfrak a}^{\bf c}+
\sum\limits_{\widetilde{\alpha}\in\widetilde{\triangle}}
\mathfrak g^{\bf c}_{\widetilde{\alpha}}$ and 
${\rm dim}_{\bf c}\mathfrak g^{\bf c}_{\widetilde{\alpha}}=1$ for each 
$\widetilde{\alpha}\in\widetilde{\triangle}$.  Also, we have 
$\triangle=\{\widetilde{\alpha}\vert_{{\mathfrak a}^{\bf c}}\,\vert\,
\widetilde{\alpha}\in\widetilde{\triangle}\}\setminus\{0\}$, 
$\mathfrak q^{\bf c}_{\alpha}=
(\sum\limits_{\widetilde{\alpha}\in\widetilde{\triangle}\,\,{\rm s.t.}\,\,
\widetilde{\alpha}\vert_{\mathfrak a^{\bf c}}=\pm\alpha})\cap
\mathfrak q^{\bf c}$ ($\alpha\in\triangle$) and $\mathfrak h^{\bf c}_{\alpha}=
(\sum\limits_{\widetilde{\alpha}\in\widetilde{\triangle}\,\,{\rm s.t.}\,\,
\widetilde{\alpha}\vert_{\mathfrak a^{\bf c}}=\pm\alpha})\cap
\mathfrak h^{\bf c}$ ($\alpha\in\triangle$).  
The following fact is well-known.  

\vspace{0.5truecm}

\noindent
{\bf Lemma 3.1.} {\sl For each $\alpha\in\triangle$, 
$\alpha(\mathfrak a\cap\mathfrak p)\subset{\bf R}$ and 
$\alpha(\mathfrak a\cap\mathfrak f)\subset\sqrt{-1}{\bf R}$.}

\vspace{0.5truecm}

\noindent
{\it Remark 3.1.} Each element of $\mathfrak a\cap\mathfrak p$ (resp. 
$\mathfrak a\cap\mathfrak f$) is called a {\it hyperbolic} (resp. 
{\it elliptic}) {\it element}.  

\vspace{0.5truecm}

Take $E_{\widetilde{\alpha}}(\not=0)\in
\mathfrak g^{\bf c}_{\widetilde{\alpha}}$ for each 
$\widetilde{\alpha}\in\widetilde{\triangle}$ and set 
$Z_{\widetilde{\alpha}}:=c_{\widetilde{\alpha}}(E_{\widetilde{\alpha}}
+\sigma E_{\widetilde{\alpha}})$ and 
$Y_{\widetilde{\alpha}}:=c_{\widetilde{\alpha}}
(E_{\widetilde{\alpha}}-\sigma E_{\widetilde{\alpha}})$, where 
$c_{\widetilde{\alpha}}$ is one of two solutions of the complex equation 
$$z^2=\frac{\alpha(a_{\alpha})}
{B^{\bf c}(E_{\widetilde{\alpha}}-\sigma E_{\widetilde{\alpha}},
E_{\widetilde{\alpha}}-\sigma E_{\widetilde{\alpha}})}.$$
Then we have ${\rm ad}(a)Z_{\widetilde{\alpha}}=\widetilde{\alpha}(a)
Y_{\widetilde{\alpha}}$ and ${\rm ad}(a)Y_{\widetilde{\alpha}}
=\widetilde{\alpha}(a)Z_{\widetilde{\alpha}}$ for any 
$a\in\mathfrak a^{\bf c}$.  Hence we have $Z_{\widetilde{\alpha}}\in
\mathfrak h^{\bf c}_{\widetilde{\alpha}\vert_{\mathfrak a^{\bf c}}}$ and 
$Y_{\widetilde{\alpha}}\in
\mathfrak q^{\bf c}_{\widetilde{\alpha}\vert_{\mathfrak a^{\bf c}}}$.  
Furthermore, for $\alpha\in\mathfrak a^{\bf c}$, it is shown that 
$\mathfrak h^{\bf c}_{\alpha}$ (resp. $\mathfrak q^{\bf c}_{\alpha}$) is 
spanned by $\{Z_{\widetilde{\alpha}}\,\vert\,\widetilde{\alpha}\in
\widetilde{\triangle}\,\,{\rm s.t.}\,\,\widetilde{\alpha}
\vert_{\mathfrak a^{\bf c}}=\alpha\}$ (resp. $\{Y_{\widetilde{\alpha}}\,\vert
\,\widetilde{\alpha}\in\widetilde{\triangle}\,\,{\rm s.t.}\,\,
\widetilde{\alpha}\vert_{\mathfrak a^{\bf c}}=\alpha\}$).  For each 
$\alpha\in\triangle$, define $a_{\alpha}\in\mathfrak a^{\bf c}$ by 
$\alpha(a)=B^{\bf c}(a_{\alpha},a)$ ($a\in\mathfrak a^{\bf c}$).  Then 
$[Z_{\widetilde{\alpha}},Y_{\widetilde{\alpha}}]=\alpha(a_{\alpha})a_{\alpha}$ 
is shown.  
L. Verhoczki [38] investigated the shape operators of orbits of the isotropy 
action of a Riemannian symmetric space of compact type.  By applying his 
method of investigation to the isotropy action of the anti-Kaehlerian 
symmetric space $G^{\bf c}/H^{\bf c}$, we prove the following fact for orbits 
of the isotropy action of $G/H$.  

\vspace{0.5truecm}

\noindent
{\bf Proposition 3.2.} {\sl Let $M$ be an orbit of the isotropy action (i.e., 
the $H$-action) on $G/H$ through $x:=\exp_G(w)H$ ($w\in\mathfrak q$ s.t. 
${\rm ad}(w)\,:\,$ semi-simple) and $A$ be the shape tensor of $M$.  
For simplicity, set $g:=\exp_G(w)$.  Let 
$\mathfrak a$ be a Cartan subspace of $\mathfrak q$ containing 
$w$ and $\mathfrak q^{\bf c}=\mathfrak a^{\bf c}
+\sum\limits_{\alpha\in\triangle_+}
\mathfrak q^{\bf c}_{\alpha}$ be the root space decomposition with respect to 
$\mathfrak a^{\bf c}$.  Then the following statements {\rm (i)} and 
{\rm (ii)} hold:

{\rm (i)} $g_{\ast}^{-1}(T_xM)^{\bf c}=\sum\limits_{\alpha\in\triangle_+\,\,
{\rm s.t.}\,\,\alpha(w)\notin\sqrt{-1}\pi{\bf Z}}\mathfrak q^{\bf c}_{\alpha}$ 
and $g_{\ast}^{-1}(T^{\perp}_xM)^{\bf c}=\mathfrak a^{\bf c}+$\newline
$\sum\limits_{\alpha\in \triangle_+\,\,{\rm s.t.}\,\,\alpha(w)\in\sqrt{-1}\pi
{\bf Z}}\mathfrak q^{\bf c}_{\alpha}$ hold.  In particular, if $M$ is a 
principal orbit, then we have $g_{\ast}^{-1}(T_xM)^{\bf c}
=\sum\limits_{\alpha\in\triangle_+}\mathfrak q_{\alpha}^{\bf c}$ and 
$g_{\ast}^{-1}(T^{\perp}_xM)^{\bf c}=\mathfrak a^{\bf c}$.  

{\rm (ii)} Let $H_x$ be the isotropy group of $H$ at $x$ and set 
$H_x(g_{\ast}\mathfrak a):=\{h_{\ast x}g_{\ast}a\,\vert\,
a\in\mathfrak a,\,h\in H_x\}$.  
Then $H_x(g_{\ast}\mathfrak a)$ is open in $T^{\perp}_xM$ and, 
for any $v:=h_{\ast x}g_{\ast}a\in H_x(g_{\ast}a)$ 
($a\in\mathfrak a,\,h\in H_x$), we have 
$A^{\bf c}_v\vert_{h_{\ast x}g_{\ast}\mathfrak q^{\bf c}_{\alpha}}
=-\frac{\sqrt{-1}\alpha(a)}
{\tan(\sqrt{-1}\alpha(w))}{\rm id}$ 
{\rm(}$\alpha\in\triangle_+$ s.t. $\alpha(w)\notin\sqrt{-1}\pi{\bf Z}${\rm)}, 
where $A^{\bf c}$ is the complexification of $A$.}

\vspace{0.5truecm}

\noindent
{\it Proof.} First we shall show the statement (i) by imitating the proof of 
Proposition 3 in [38].  Let $M^{\bf c}$ be the extrinsic complexification of 
$M$, that is, $M^{\bf c}:=H^{\bf c}\cdot x\,(\subset G^{\bf c}/H^{\bf c})$, 
where $G/H$ is identified with $G(eH^{\bf c})$.  We shall investigate 
$T_x(M^{\bf c})$ instead of $(T_xM)^{\bf c}$ because $(T_xM)^{\bf c}$ is 
identified with $T_x(M^{\bf c})$.  Let $a_{\alpha}$ ($\alpha\in\triangle$), 
$\widetilde{\triangle},\,Z_{\widetilde{\alpha}}$ and $Y_{\widetilde{\alpha}}$ 
($\widetilde{\alpha}\in\widetilde{\triangle}$) be the above quantities defined 
for $\mathfrak a$ and a Cartan subalgebra $\widetilde{\mathfrak a}$ of 
$\mathfrak g$ containing $\mathfrak a$.  Let $\widetilde{\alpha}\in
\widetilde{\triangle}$ and $\alpha:=\widetilde{\alpha}
\vert_{{\mathfrak a}^{\bf c}}$.  Since $[Z_{\widetilde{\alpha}},w]
=-\alpha(w)Y_{\widetilde{\alpha}}$ and 
$[Z_{\widetilde{\alpha}},Y_{\widetilde{\alpha}}]=\alpha(a_{\alpha})
a_{\alpha}$, we have 
$$\frac{d}{dt}\vert_{t=0}{\rm Ad}_{G^{\bf c}}(\exp tZ_{\widetilde{\alpha}})w
=-\alpha(w)Y_{\widetilde{\alpha}},$$
where ${\rm Ad}_{G^{\bf c}}$ is the adjoint representation of $G^{\bf c}$.  
Hence we have 
$$T_w{\rm Ad}_{G^{\bf c}}(H^{\bf c})w
=\sum_{\alpha\in\triangle_+\,\,{\rm s.t.}\,\,\alpha(w)\not=0}
\mathfrak q^{\bf c}_{\alpha}.$$
Denote by ${\rm Exp}$ the exponential map of 
the anti-Kaehlerian symmetric space $(G^{\bf c}/H^{\bf c},J,\langle\,\,,\,\,
\rangle)$.  Assume that $\alpha(w)\not=0$.  Define a complex geodesic 
variation $\delta:{\bf C}^2\to G^{\bf c}/H^{\bf c}$ of the complex geodesic 
$\gamma^{\bf c}_w\,(z\mapsto{\rm Exp}(zw))$ by 
$$\delta(z,u):={\rm Exp}\left(z(\cos u\cdot w+\sin u
\sqrt{\frac{\langle w,w\rangle}{\langle Y_{\widetilde{\alpha}},
Y_{\widetilde{\alpha}}\rangle}}Y_{\widetilde{\alpha}})\right)$$
($(z,u)\in{\bf C}^2$).  Set $W:=\frac{\partial\delta}{\partial u}\vert_{u=0}$, 
which is a complex Jacobi field along $\gamma_w^{\bf c}$.  Hence it follows 
from $(2.2)$ that 
$$W_1=\frac{\sin(\sqrt{-1}\alpha(w))}{\sqrt{-1}\alpha(w)}
\sqrt{\frac{\langle w,w\rangle}{\langle Y_{\widetilde{\alpha}},
Y_{\widetilde{\alpha}}\rangle}}g_{\ast}Y_{\widetilde{\alpha}}.$$
On the other hand, we have 
$W_1=(d{\rm Exp})_w(\sqrt{\frac{\langle w,w\rangle}
{\langle Y_{\widetilde{\alpha}},Y_{\widetilde{\alpha}}\rangle}}
Y_{\widetilde{\alpha}})$.  Hence we have 
$$(d{\rm Exp})_w(Y_{\widetilde{\alpha}})=
\frac{\sin(\sqrt{-1}\alpha(w))}{\sqrt{-1}\alpha(w)}g_{\ast}
Y_{\widetilde{\alpha}}. \leqno{(3.1)}$$
Since $M^{\bf c}={\rm Exp}({\rm Ad}_{G^{\bf c}}(H^{\bf c})w)$, we have 
$T_x(M^{\bf c})=(d{\rm Exp})_w(T_w({\rm Ad}_{G^{\bf c}}(H^{\bf c})w))$.  
Hence the relations in the statement (i) follow from $(3.1)$.  

Next we shall show the statement (ii).  
The $H_x$-action on $T_x(G/H)$ preserves $T_xM$ and $T_x^{\perp}M$ 
invariantly, respectively.  
The $H_x$-action on $T^{\perp}_xM$ is so-called 
slice representation and it is equivalent to an $s$-representation 
(the linear isotropy representation of a pseudo-Riemannian symmetric space).  
Therefore $H_x(g_{\ast}\mathfrak a)$ is open in $T^{\perp}_xM$.  
In the sequel, we shall show the remaining part of the statement (ii) 
by imitating 
the proof of Theorem 1 in [38] for the isotropy action of a Riemannian 
symmetric space of compact type.  
Denote by $\widehat A$ the shape tensor of $M^{\bf c}$.  Under the 
identification of $(T_xM)^{\bf c}$ with $T_x(M^{\bf c})$, the complexified 
shape operator $A^{\bf c}_w$ is identified with $\widehat A_w$.  Hence we 
suffice to investigate $\widehat A_w$ instead of $A^{\bf c}_w$.  
Let $\alpha$ be an element of $\triangle_+$ with 
$\alpha(w)\notin\sqrt{-1}\pi{\bf Z}$.  Take $\widetilde{\alpha}_1\in
\widetilde{\triangle}$ with $\widetilde{\alpha}_1\vert_{\mathfrak a^{\bf c}}
=\alpha$.  Also, in case of $2\alpha\in\triangle$, $\widetilde{\alpha}_2\in
\widetilde{\triangle}$ with $\widetilde{\alpha}_2\vert_{\mathfrak a^{\bf c}}
=2\alpha$.  Set $\widehat{\mathfrak h}^{\bf c}_{\alpha}
:=\mathfrak z_{\mathfrak h^{\bf c}}(\mathfrak a^{\bf c})
+\mathfrak h^{\bf c}_{\alpha}+\mathfrak h^{\bf c}_{2\alpha}$ 
($\mathfrak h^{\bf c}_{2\alpha}=\{0\}$ in case of $2\alpha\notin\triangle$) 
and $\widehat H^{\bf c}_{\alpha}:=\exp_{G^{\bf c}}
(\widehat{\mathfrak h}^{\bf c}_{\alpha})$.  Easily we can show 
$${\rm Ad}_{G^{\bf c}}(\exp\,zZ_{\widetilde{\alpha}_k})a_{\alpha}
=\cos(k^2z\alpha(a_{\alpha}))a_{\alpha}-\frac1k
\sin(k^2z\alpha(a_{\alpha}))Y_{\widetilde{\alpha}_k}\quad(k=1,2).$$
From this relation, it follows that ${\rm Ad}(\widehat H^{\bf c}_{\alpha})
(a_{\alpha})$ is a complex hypersurface in 
$\widehat{\mathfrak q}^{\bf c}_{\alpha}:={\bf C}a_{\alpha}
+\mathfrak q^{\bf c}_{\alpha}+\mathfrak q^{\bf c}_{2\alpha}$ 
($\mathfrak q^{\bf c}_{2\alpha}=\{0\}$ in case of $2\alpha\notin\triangle$).  
On the other hand, it is clear that ${\rm Ad}(\widehat H^{\bf c}_{\alpha})
(a_{\alpha})$ is contained in the complex hypersphere 
$(B^{\bf c}\vert_{\mathfrak q^{\bf c}_{\alpha}
\times\mathfrak q^{\bf c}_{\alpha}})({\bf z},{\bf z})
=B^{\bf c}(a_{\alpha},a_{\alpha})$ of 
$\widehat{\mathfrak q}^{\bf c}_{\alpha}$.  Hence 
${\rm Ad}(\widehat H^{\bf c}_{\alpha})(a_{\alpha})$ coincides with this 
complex hypersphere.  The vector $w$ is expressed as 
$w=\frac{\alpha(w)}{\alpha(a_{\alpha})}a_{\alpha}+b$ for some 
$b\in\alpha^{-1}(0)$.  Then we have 
$${\rm Ad}_{G^{\bf c}}(\exp zZ_{\widetilde{\alpha}_k})w=b+
\frac{\alpha(w)}{\alpha(a_{\alpha})}(\cos(k^2z\alpha(a_{\alpha}))a_{\alpha}
-\frac1k\sin(k^2z\alpha(a_{\alpha}))Y_{\widetilde{\alpha}_k})$$
$(k=1,2)$.  From this relation, it follows that 
${\rm Ad}(\widehat H^{\bf c}_{\alpha})(w)$ coincides with the complex 
hypersphere 
$(B^{\bf c}\vert_{{\widehat{\mathfrak q}}^{\bf c}_{\alpha}
\times{\widehat{\mathfrak q}}^{\bf c}_{\alpha}})({\bf z}-b,{\bf z}-b)
=\frac{\alpha(w)^2}{\alpha(a_{\alpha})}$ of $b+
{\widehat{\mathfrak q}}^{\bf c}_{\alpha}$.  
Set ${\widehat Q}^{\bf c}_{\alpha}:={\rm Exp}
({\widehat{\mathfrak q}}^{\bf c}_{\alpha})$ and 
${\widehat Q}^{\bf c}_{\alpha}(b):={\rm Exp}
(b+{\widehat{\mathfrak q}}^{\bf c}_{\alpha})$.  
It is easy to show that ${\widehat Q}^{\bf c}_{\alpha}$ is a totally geodesic 
complex rank one anti-Kaehlerian symmetric space in $G^{\bf c}/H^{\bf c}$.  
Furthermore, by imitating the proof of Proposition 4 in [38], it is shown that 
${\widehat Q}^{\bf c}_{\alpha}(b)$ is a totally geodesic complex rank one 
anti-Kaehlerian symmetric space and it is isometric to 
${\widehat Q}^{\bf c}_{\alpha}$.  In fact, a map 
$\phi:{\widehat Q}^{{\bf c}'}_{\alpha}\to{\widehat Q}^{\bf c'}_{\alpha}(b)$ 
defined by $\phi({\rm Exp}{\bf z})={\rm Exp}({\bf z}+b)$ (${\bf z}\in
{\widehat Q}^{{\bf c}'}_{\alpha}$) is an isometry.  
Since ${\rm Ad}({\widehat H}^{\bf c}_{\alpha})(w)$ is equal to the complex 
hypersphere of complex radius $\sqrt{\frac{\alpha(w)^2}{\alpha(a_{\alpha})}}$ 
of $b+{\widehat{\mathfrak q}}^{\bf c}_{\alpha}$, $\widehat H^{\bf c}_{\alpha}
\cdot x$ is a complex geodesic hypersphere of complex radius 
$\sqrt{\alpha(w)}$ in $\widehat Q^{\bf c}_{\alpha}(b)$.  Set 
${\widehat Q}^{{\bf c}'}_{\alpha}:={\rm Exp}({\mathfrak a}^{\bf c}
+{\mathfrak q}^{\bf c}_{\alpha}+{\mathfrak q}^{\bf c}_{2\alpha})$, which is 
isometric to the anti-Kaehlerian product $\widehat Q^{\bf c}_{\alpha}(b)\times
{\bf C}^{r-1}$ ($r:={\rm rank}(G/H)$).  

We have ${\widehat H}^{\bf c}_{\alpha}\cdot x\subset M^{\bf c}\cap
{\widehat Q}^{\bf c}_{\alpha}(b)\subset M^{\bf c}\cap
{\widehat Q}^{{\bf c}'}_{\alpha}$.  Also, since 
$T_x(M^{\bf c})=g_{\ast}(\sum\limits_{\alpha\in\triangle_+\,\,{\rm s.t.}\,\,
\alpha(w)\notin\sqrt{-1}\pi{\bf Z}}{\mathfrak q}^{\bf c}_{\alpha})$ and 
$T_x{\widehat Q}^{{\bf c}'}_{\alpha}=g_{\ast}(\mathfrak a^{\bf c}
+\mathfrak q^{\bf c}_{\alpha}+\mathfrak q^{\bf c}_{2\alpha})$, we have 
$T_x(M^{\bf c}\cap{\widehat Q}^{{\bf c}'}_{\alpha})
=\mathfrak q^{\bf c}_{\alpha}+\mathfrak q^{\bf c}_{2\alpha}$ and hence 
${\rm dim}\,T_x(M^{\bf c}\cap
{\widehat Q}^{{\bf c}'}_{\alpha})={\rm dim}\,({\widehat H}^{\bf c}_{\alpha}
\cdot x)$.  
Therefore ${\widehat H}^{\bf c}_{\alpha}\cdot x$ is a component of 
$M^{\bf c}\cap{\widehat Q}^{{\bf c}'}_{\alpha}$.  Denote by $\overline A$ 
the shape tensor of ${\widehat H}^{\bf c}_{\alpha}\cdot x\hookrightarrow 
{\widehat Q}^{{\bf c}'}_{\alpha}$.  Since 
${\widehat Q}^{{\bf c}'}_{\alpha}$ is totally geodesic in 
$G^{\bf c}/H^{\bf c}$ and 
$T^{\perp}_x(M^{\bf c})$ contains the normal space of 
${\widehat H}^{\bf c}_{\alpha}\cdot x$ in 
${\widehat Q}^{{\bf c}'}_{\alpha}$, it follows 
from pseudo-Riemannian version of Lemma 6 of [38] that 
$\widehat A_{g_{\ast}a_{\alpha}}$ preserves 
$T_x({\widehat H}^{\bf c}_{\alpha}\cdot x)$ invariantly and that 
$\widehat A_{g_{\ast}a_{\alpha}}=\overline A_{g_{\ast}a_{\alpha}}$ on 
$T_x({\widehat H}^{\bf c}_{\alpha}\cdot x)$.  
Let $\phi$ be the above isometry of $\widehat Q^{\bf c}_{\alpha}$ onto 
$\widehat Q^{\bf c}_{\alpha}(b)$.  Set $r_0:=\frac{\alpha(w)}
{\alpha(a_{\alpha})}$ and denote by ${\overline A}'$ the shape tensor of 
${\widehat H}^{\bf c}_{\alpha}\cdot(r_0a_{\alpha})\hookrightarrow
{\widehat Q}^{{\bf c}'}_{\alpha}$.  Clearly we have 
$\phi({\widehat H}^{\bf c}_{\alpha}\cdot(r_0a_{\alpha}))
={\widehat H}^{\bf c}_{\alpha}\cdot x$ and $\phi_{\ast}
((\exp_{G^{\bf c}}(r_0a_{\alpha}))_{\ast}(a_{\alpha}))=g_{\ast}a_{\alpha}$.  
Hence we have $\overline A_{g_{\ast}a_{\alpha}}=\phi_{\ast}\circ
{\overline A}'_{(\exp_{G^{\bf c}}(r_0a_{\alpha}))_{\ast}(a_{\alpha})}
\circ\phi_{\ast}^{-1}$.  For simplicity, set $\overline g:=\exp_{G^{\bf c}}
(r_0a_{\alpha})$.  Now we shall investigate 
${\overline A}'_{\overline g_{\ast}a_{\alpha}}$.  
Define 
a complex geodesic variation $\delta:{\bf C}^2\to G^{\bf c}/H^{\bf c}$ by 
$$\delta(z,u):={\rm Exp}(z(r_0\cos u\cdot a_{\alpha}+
\sqrt{\frac{r_0^2\langle a_{\alpha},a_{\alpha}\rangle}
{\langle Y_{\widetilde{\alpha}_1},Y_{\widetilde{\alpha}_1}\rangle}}\sin u\cdot 
Y_{\widetilde{\alpha}_1}))\,\,\,\,((z,u)\in{\bf C}^2).$$
Set $W:=\frac{\partial\delta}{\partial u}\vert_{u=0}$.  Since $W$ is a complex 
Jacobi field along $\gamma_{r_0a_{\alpha}}^{\bf c}$, 
it follows from $(2.2)$ that 
$$W_z=\frac{\sin(\sqrt{-1}z\alpha(r_0a_{\alpha}))}
{\sqrt{-1}\alpha(r_0a_{\alpha})}\sqrt{\frac{r_0^2\langle a_{\alpha},
a_{\alpha}\rangle}{\langle Y_{\widetilde{\alpha}_1},
Y_{\widetilde{\alpha}_1}\rangle}}(P_{\gamma^{\bf c}_{r_0a_{\alpha}}})_{0,z}
(Y_{\widetilde{\alpha}_1}).
\leqno{(3.2)}
$$
We have 
$$\begin{array}{l}
\displaystyle{
\widetilde{\nabla}_{\frac{\partial\delta}{\partial u}\vert_{z=1,u=0}}
\frac{\partial\delta}{\partial z}=\widetilde{\nabla}_{\frac{\partial\delta}
{\partial z}\vert_{z=1,u=0}}\frac{\partial\delta}{\partial u}=W'_1}\\
\displaystyle{
=\cos(\sqrt{-1}\alpha(r_0a_{\alpha}))\sqrt{\frac{r_0^2\langle a_{\alpha},
a_{\alpha}\rangle}{\langle Y_{\widetilde{\alpha}_1},
Y_{\widetilde{\alpha}_1}\rangle}}\overline g_{\ast}Y_{\widetilde{\alpha}_1}
\in T_{{\rm Exp}(r_0a_{\alpha})}
\widehat H^{\bf c}_{\alpha}\cdot(r_0a_{\alpha})}
\end{array}$$
and hence 
$${\overline A}'_{\overline g_{\ast}(r_0a_{\alpha})}W_1
=-\cos(\sqrt{-1}\alpha(r_0a_{\alpha}))
\sqrt{\frac{r_0^2\langle a_{\alpha},a_{\alpha}\rangle}
{\langle Y_{\widetilde{\alpha}_1},Y_{\widetilde{\alpha}_1}\rangle}}
\overline g_{\ast}Y_{\widetilde{\alpha}_1},$$
which together with $(3.2)$ and $\alpha(b)=0$ deduces 
$${\overline A}'_{\overline g_{\ast}a_{\alpha}}
\overline g_{\ast}Y_{\widetilde{\alpha}_1}
=-\frac{\sqrt{-1}\alpha(a_{\alpha})}{\tan(\sqrt{-1}\alpha(w))}
\overline g_{\ast}Y_{\widetilde{\alpha}_1}.$$
Therefore we have 
$$\widehat A_{g_{\ast}a_{\alpha}}g_{\ast}Y_{\widetilde{\alpha}_1}
=-\frac{\sqrt{-1}\alpha(a_{\alpha})}{\tan(\sqrt{-1}\alpha(w))}
g_{\ast}Y_{\widetilde{\alpha}_1}.$$
Similarly we have 
$$\widehat A_{g_{\ast}a_{\alpha}}g_{\ast}Y_{\widetilde{\alpha}_2}
=-\frac{2\sqrt{-1}\alpha(a_{\alpha})}{\tan(2\sqrt{-1}\alpha(w))}
g_{\ast}Y_{\widetilde{\alpha}_2}.$$
Take $\bar b\in\alpha^{-1}(0)$.  
Since $\widehat Q^{\bf c}_{\alpha}(b)$ is totally geodesic and 
$T^{\perp}\widehat Q^{\bf c}_{\alpha}(b)
\vert_{\widehat H^{\bf c}_{\alpha}\cdot x}\cap T^{\perp}M
\vert_{\widehat H^{\bf c}_{\alpha}\cdot x}$ is parallel along 
$\widehat H^{\bf c}_{\alpha}\cdot x$ with respect to the normal connection of 
$\widehat Q^{\bf c}_{\alpha}(b)\hookrightarrow G^{\bf c}/H^{\bf c}$, we have 
$$\widehat A_{g_{\ast}\bar b}g_{\ast}Y_{\widetilde{\alpha}_1}
=\widehat A_{g_{\ast}\bar b}g_{\ast}Y_{\widetilde{\alpha}_2}=0.$$
Take an arbitrary $a\in\mathfrak a$.  We can express as 
$a=\frac{\alpha(a)}{\alpha(a_{\alpha})}a_{\alpha}+\widehat b$ 
for some $\widehat b\in\alpha^{-1}(0)$.  
Thus, for each $a\in\mathfrak a$, we have 
$$\widehat A_{g_{\ast}a}\vert_{g_{\ast}\mathfrak q^{\bf c}_{\beta}}
=-\frac{\sqrt{-1}\beta(a)}{\tan(\sqrt{-1}\beta(w))}{\rm id}
\,\,\,(\beta\in\triangle_+\,\,{\rm s.t.}\,\,\beta(w)\notin\sqrt{-1}
\pi{\bf Z}). 
$$
Take an arbitrary $h_{\ast x}g_{\ast}a\in H_x(g_{\ast}\mathfrak a)$ 
($a\in\mathfrak a,\,h\in H_x$).  Since $h$ is an 
isometry of $G^{\bf c}/H^{\bf c}$, we have 
$\widehat A_{h_{\ast x}g_{\ast}a}=h_{\ast x}\circ\widehat A_{g_{\ast}a}
\circ h_{\ast x}^{-1}$.  
Hence we have 
$$\widehat A_{h_{\ast x}g_{\ast}a}\vert_{h_{\ast x}g_{\ast}
\mathfrak q^{\bf c}_{\beta}}
=-\frac{\sqrt{-1}\beta(a)}{\tan(\sqrt{-1}\beta(w))}{\rm id}
\,\,\,(\beta\in\triangle_+\,\,{\rm s.t.}\,\,\beta(w)\notin\sqrt{-1}
\pi{\bf Z}).$$
Therefore, we obtain the relation in the statement (ii).  
\hspace{1.5truecm}q.e.d.

\section{Shape operators of partial tubes}
In this section, we investigate the shape operators of partial tubes over 
a pseudo-Riemannian submanifold with section in a (semi-simple) 
pseudo-Riemannian symmetric space $G/H$ 
equipped with the metric induced from the Killing form of 
$\mathfrak g:={\rm Lie}\,G$.  
Let $M$ be a pseudo-Riemannian submanifold with section in $G/H$, that is, 
for each $x=gH$ of $M$, $g_{\ast}^{-1}T^{\perp}_xM$ is a Lie triple system.  
Let $t(M)$ be a connected submanifold in the normal bundle $T^{\perp}M$ of 
$M$ such that, for any curve $c:[0,1]\to M$, $P^{\perp}_c(t(M)\cap 
T^{\perp}_{c(0)}M)=t(M)\cap T^{\perp}_{c(1)}M$ holds, where 
$P^{\perp}_c$ is the parallel transport along $c$ with respect to 
the normal connection.  Denote by $F$ the set of all critical points of the 
normal exponential map $\exp^{\perp}$ of $M$.  Assume that $t(M)\cap F
=\emptyset$.  Then the restriction $\exp^{\perp}\vert_{t(M)}$ of 
$\exp^{\perp}$ to $t(M)$ is an immersion of $t(M)$ into $G/H$.  
Assume that $\exp^{\perp}\vert_{t(M)}\,:\,t(M)
\hookrightarrow G/H$ is a pseudo-Riemannian submanifold.  
Then we call $t(M)$ a {\it partial tube over} $M$.  
Define a distribution $D^V$ on $t(M)$ by 
$D^V_v=T_v(t(M)\cap T^{\perp}_{\pi(v)}M)$ ($v\in t(M)$), where 
$\pi$ is the bundle projection of $T^{\perp}M$.  We call this distribution 
a {\it vertical distribution} on $t(M)$.  Let $X\in T_{\pi(v)}M$.  Take 
a curve $c$ in $M$ with $\dot c(0)=X$.  
Let $\widetilde v$ be a parallel normal vector field along $c$ with 
$\widetilde v(0)=v$.  
Denote by $\widetilde X_v$ the velocity vector $\dot{\widetilde v}(0)$ of 
the curve $\widetilde v$ in $T^{\perp}M$ at $0$.  
We call $\widetilde X_v$ the {\it horizontal lift of} $X$ {\it to} 
$v$.  Define a distribution $D^H$ on $t(M)$ by $D^H_v=\{\widetilde X_v
\,\vert\,X\in T_{\pi(v)}M\}$ ($v\in t(M)$).  We call this distribution 
a {\it horizontal distribution} on $t(M)$.  
From $(2.1)$, we have 
$$\exp^{\perp}_{\ast}(\widetilde X_v)=P_{\gamma_v}\left(D^{co}_vX-D^{si}_v
(A_vX)\right). \leqno{(4.1)}$$
Assume that $t(M)$ is contained in the 
$\varepsilon$-tube $t_{\varepsilon}(M):=\{v\in T^{\perp}M\,\vert\,
\frac{\langle v,v\rangle}{\sqrt{\vert\langle v,v\rangle\vert}}=\varepsilon\}$ 
($\varepsilon\not=0$).  Define a subbundle $D^{\perp}$ of the normal bundle 
$T^{\perp}t(M)$ of $t(M)$ by $D^{\perp}_v:=T^{\perp}_vt(M)\cap 
T_v(t_{\varepsilon}(M))$ ($v\in t(M)$).  
Clearly we have $T_vt(M)=D^H_v\oplus D^V_v$ (orthogonal direct sum) 
and $T^{\perp}_vt(M)=D^{\perp}_v\oplus{\rm Span}\{\dot{\bar{\gamma}}_v(1)\}$ 
(orthogonal direct sum), where $\bar{\gamma}_v$ is defined by 
$\bar{\gamma}_v(t):=tv$.  Denote by $A$ (resp. $A^t$) the shape tensor of $M$ 
(resp. $t(M)$).  Also, denote by $A^x$ that of a submanifold $t(M)\cap 
T^{\perp}_xM$ in $\exp^{\perp}(T^{\perp}_xM)$ immersed by $\exp^{\perp}
\vert_{t(M)\cap T^{\perp}_xM}$.  In the sequel, we omit 
$\exp^{\perp}_{\ast}$.  For a real analytic 
function $F$ and $v\in T_{gH}(G/H)$, we denote the operator 
$g_{\ast}\circ F({\rm ad}(g_{\ast}^{-1}v))\circ g_{\ast}^{-1}$ by 
$F({\rm ad}(v))$ for simplicity.  
Then, by imitating the proof of Proposition 3.1 in [19], we can show 
the the following relations.  

\vspace{0.5truecm}

\noindent
{\bf Proposition 4.1.} 
{\sl Let $v\in t(M)$ and $w\in D^{\perp}_v$.  Also, let $\pi(v)=g_1H,\,
g_2:=\exp_G(g_{1\ast}^{-1}v)$ and $g:=g_1g_2g_1^{-1}$, where 
$\exp_G$ is the exponential map of the Lie group $G$.  

{\rm(}i{\rm)} For $Y\in D^V_v$, we have 
$$A^t_{g_{\ast}v}Y=A^{\pi(v)}_{g_{\ast}v}Y,\quad 
A^t_wY=A^{\pi(v)}_wY. \leqno{(4.2)}$$

{\rm(}ii{\rm)} Assume that 
${\rm Span}\{g_{1\ast}^{-1}v,(g_1g_2)_{\ast}^{-1}w\}$ is abelian.  
Then, for $X\in T_{\pi(v)}M$, we have 
$$\begin{array}{l}
\displaystyle{A^t_w\widetilde X_v=
\sqrt{-1}{\rm ad}(g_{\ast}^{-1}w)\sin(\sqrt{-1}{\rm ad}(v))(X)}\\
\hspace{1.6truecm}\displaystyle{
-\frac{\sqrt{-1}\sin(\sqrt{-1}{\rm ad}(v))}{{\rm ad}(v)}
(A_{g_{\ast}^{-1}w}X)}\\
\hspace{1.6truecm}\displaystyle{+\left(
\frac{\cos(\sqrt{-1}{\rm ad}(v))-{\rm id}}{{\rm ad}(v)}+
\frac{\sqrt{-1}\sin(\sqrt{-1}{\rm ad}(v))+{\rm ad}(v)}{{\rm ad}(v)^2}\right)}\\
\hspace{6.5truecm}\displaystyle{\times{\rm ad}(g_{\ast}^{-1}w)
(A_vX).}
\end{array}\leqno{(4.3)}$$
}

\vspace{0.5truecm}

\noindent
{\it Remark 4.1.} The parallel translation $P_{\gamma_v}$ along $\gamma_v$ 
is equal to $g_{\ast}$.  

\section{Proper complex equifocality} 
In this section, we investigate the proper complex equifocality of a complex 
equifocal submanifold in a pseudo-Riemannian symmetric space.  
Let $G/H$ be a (semi-simple) pseudo-Riemannian symmetric space and $R$ be 
the curvature tensor of $G/H$.  First we prepare the following lemma for 
a curvature-adapted submanifold with flat section such that the normal 
holonomy group is trivial.  

\vspace{0.5truecm}

\noindent
{\bf Lemma 5.1.} {\sl Let $M$ be a curvature-adapted submanifold in $G/H$ 
with flat section such that the normal holonomy group is trivial.  
Assume that, for any normal vector 
$v$ of $M$, $A_v$ and ${\rm ad}(g_{\ast}^{-1}v)$ are semi-simple, where $A$ 
is the shape tensor of $M$ and $g$ is an element of $G$ 
such that $gH$ is the base point of $v$.  Then, for any $x\in M$, 
$\{A_v\,\vert\,v\in T_x^{\perp}M\}\cup\{R(\cdot,v)v\vert_{T_xM}\,\vert\,
v\in T^{\perp}_xM\}$ is a commuting family of linear transformations of 
$T_xM$.}

\vspace{0.5truecm}

\noindent
{\it Proof.} Let $v_i\in T^{\perp}_xM$ ($i=1,2$).  Since $M$ has flat section, 
$R(\cdot,v_1)v_1\vert_{T_xM}$ and $R(\cdot,v_2)v_2\vert_{T_xM}$ 
commute with each other.  Since $M$ has flat section and the normal holonomy 
group is trivial, $A_{v_1}$ and $A_{v_2}$ commute with each other.  
In the sequel, we shall show that $R(\cdot,v_1)v_1\vert_{T_xM}$ and $A_{v_2}$ 
commute with each other.  
Let $x=gH$.  Since $g_{\ast}^{-1}T^{\perp}_xM$ is abelian and, for any 
$v\in T^{\perp}_xM$, 
${\rm ad}(g_{\ast}^{-1}v)$ is semi-simple, there exists a Cartan subspace 
$\mathfrak a$ of $\mathfrak q(=T_{eH}(G/H))$ containing 
$\mathfrak b:=g_{\ast}^{-1}(T^{\perp}_xM)$.  Let $\triangle$ be the root 
system with respect to $\mathfrak a^{\bf c}$ and set 
$\overline{\triangle}:=\{\alpha\vert_{\mathfrak b^{\bf c}}\,\vert\,
\alpha\in\triangle\,\,{\rm s.t.}\,\,\alpha\vert_{\mathfrak b^{\bf c}}
\not=0\}$.  For each $\beta\in\overline{\triangle}$, we set 
$\mathfrak q^{\bf c}_{\beta}:=\{X\in\mathfrak q^{\bf c}\,\vert\,
{\rm ad}(b)^2(X)=\beta(b)^2X\,(\forall b\in\mathfrak b^{\bf c})\}$.  Then 
we have $\mathfrak q^{\bf c}=\mathfrak z_{\mathfrak q^{\bf c}}
(\mathfrak b^{\bf c})+\sum_{\beta\in\overline{\triangle}_+}
\mathfrak q^{\bf c}_{\beta}$, where $\overline{\triangle}_+$ is the positive 
root system under some lexicographical ordering and 
$\mathfrak z_{\mathfrak q^{\bf c}}(\mathfrak b^{\bf c})$ is the centralizer 
of $\mathfrak b^{\bf c}$ in $\mathfrak q^{\bf c}$.  Consider 
$$D:=\{v\in(T^{\perp}_xM)^{\bf c}\,\vert\,{\rm Span}\{g_{\ast}^{-1}v\}
\cap\left(
\mathop{\cup}_{(\beta_1,\beta_2)\in\overline{\triangle}_+\times
\overline{\triangle}_+\,{\rm s.t.}\,\beta_1\not=\beta_2}
({\it l}_{\beta_1}\cap{\it l}_{\beta_2})\right)=\emptyset\},$$
where ${\it l}_{\beta_i}:=\beta_i^{-1}(1)$ ($i=1,2$).  It is clear that $D$ 
is open and dense in $(T^{\perp}_xM)^{\bf c}$.  Take $v\in D$.  Since 
$\beta(g_{\ast}^{-1}v)$'s ($\beta\in\overline{\triangle}_+$) are mutually 
distinct, the decomposition 
$(T_xM)^{\bf c}=g_{\ast}(\mathfrak z_{\mathfrak q^{\bf c}}(\mathfrak b^{\bf c})
\ominus\mathfrak b^{\bf c})+\sum\limits_{\beta\in\overline{\triangle}_+}
g_{\ast}\mathfrak q^{\bf c}_{\beta}$ is the eigenspace decomposition of 
$R^{\bf c}(\cdot,v)v\vert_{(T_xM)^{\bf c}}$.  Since $M$ is curvature-adapted 
and hence $[R^{\bf c}(\cdot,v)v\vert_{(T_xM)^{\bf c}},A_v^{\bf c}]=0$, we have 
$$\begin{array}{l}
\displaystyle{(T_xM)^{\bf c}=\sum_{\lambda\in{\rm Spec}\,A_v^{\bf c}}
\left(g_{\ast}(\mathfrak z_{\mathfrak q^{\bf c}}(\mathfrak b^{\bf c})\ominus
\mathfrak b^{\bf c})\cap{\rm Ker}(A_v^{\bf c}-\lambda\,{\rm id}))\right.}\\
\hspace{2.5truecm}\displaystyle{
+\sum\limits_{\lambda\in{\rm Spec}\,A_v^{\bf c}}
\sum\limits_{\beta\in\overline{\triangle}_+}
(g_{\ast}\mathfrak q^{\bf c}_{\beta}\cap{\rm Ker}(A_v^{\bf c}-\lambda\,
{\rm id})).}
\end{array}\leqno{(5.1)}$$
Suppose that $(5.1)$ does not hold for some $v_0\in(T^{\perp}_xM)^{\bf c}
\setminus D$.  Then it is easy to show that there exists a neighborhood $U$ of 
$v_0$ in $(T^{\perp}_xM)^{\bf c}$  such that $(5.1)$ does not hold for any 
$v\in U$.  Clearly we have $U\cap D=\emptyset$.  This contradicts the fact 
that $D$ is dense in $(T^{\perp}_xM)^{\bf c}$.  
Hence $(5.1)$ holds for any $v\in(T^{\perp}_xM)^{\bf c}\setminus D$.  
Therefore, $(5.1)$ holds for any $v\in(T^{\perp}_xM)^{\bf c}$.  In particular, 
$(5.1)$ holds for $v_2$.  On the other hand, the decomposition 
$(T_xM)^{\bf c}=g_{\ast}(\mathfrak z_{\mathfrak q^{\bf c}}
(\mathfrak b^{\bf c})\ominus\mathfrak b^{\bf c})+
\sum\limits_{\beta\in\overline{\triangle}_+}g_{\ast}\mathfrak q_{\beta}^{\bf c}$ is the common eigenspace decomposition of $R^{\bf c}(\cdot,v)v
\vert_{(T_xM)^{\bf c}}$'s ($v\in(T^{\perp}_xM)^{\bf c}$).  From these facts, 
we have 
$$\begin{array}{l}
\displaystyle{(T_xM)^{\bf c}=\sum_{\lambda\in{\rm Spec}\,A_{v_2}^{\bf c}}
\sum_{\mu\in{\rm Spec}\,R^{\bf c}(\cdot,v_1)v_1\vert_{(T_xM)^{\bf c}}}}\\
\hspace{2.2truecm}\displaystyle{\left(
{\rm Ker}(R^{\bf c}(\cdot,v_1)v_1\vert_{(T_xM)^{\bf c}}-\mu\,{\rm id})\cap
{\rm Ker}(A_{v_2}^{\bf c}-\lambda\,{\rm id})\right),}
\end{array}$$
which implies that $R^{\bf c}(\cdot,v_1)v_1\vert_{(T_xM)^{\bf c}}$ and 
$A_{v_2}^{\bf c}$ commute with each other.  This completes the proof.  
\hspace{1.5truecm}q.e.d.

\vspace{0.5truecm}

By this lemma, Lemma 5.3, Propositions 5.6 and 5.7 of [17] (these lemmas are 
valid even if the ambient space is a pseudo-Riemannian symmetric space), we 
can show the following fact.  

\vspace{0.5 truecm}

\noindent
{\bf Proposition 5.2.} {\sl Let $M$ be a curvature-adapted complex equifocal 
submanifold in $G/H$.  Assume that, for any normal vector $v$ of $M$, 
$A_v$ and ${\rm ad}(g_{\ast}^{-1}v)$ are semi-simple and 
that $\pm\beta(g_{\ast}^{-1}v)\notin{\rm Spec}\,A_v^{\bf c}
\vert_{g_{\ast}\mathfrak q_{\beta}^{\bf c}}$ ($\beta\in\overline{\triangle}_+$), where $g$ is an element of $G$ such that $gH$ is the base point of $v$.  
Then $M$ is proper complex equifocal.}  

\vspace{0.5 truecm}

\noindent
{\it Proof.} Let $\widetilde M:=(\pi\circ\phi)^{-1}(M)$ and denote by 
$\widetilde A$ the shape tensor of $\widetilde M$.  Fix $u\in\widetilde M$ and 
$v\in T^{\perp}_u\widetilde M$.  For simplicity, set $x(=gH)=(\pi\circ\phi)(u)$ and $\overline v:=(\pi\circ\phi)_{\ast}(v)$.  According to Lemma 5.1, it 
follows from the assumptions that $A^{\bf c}_{\overline v}$ commutes with 
$R^{\bf c}(\cdot,w)w\vert_{(T_xM)^{\bf c}}$'s ($w\in(T^{\perp}_xM)^{\bf c}$).  
Also, it follows from 
the assumptions that $A^{\bf c}_{\overline v}$ and 
$R^{\bf c}(\cdot,w)w\vert_{T_xM}$'s 
($w\in(T^{\perp}_xM)^{\bf c}$) are diagonalizable.  Hence they are 
simultaneously diagonalizable, that is, 
we have 
$$(T_x^{\perp}M)^{\bf c}=\sum_{\lambda\in{\rm Spec}\,A_{\overline v}^{\bf c}}
\sum_{\beta\in\overline{\triangle}_+}(g_{\ast}\mathfrak q_{\beta}^{\bf c}
\cap{\rm Ker}(A_{\overline v}^{\bf c}-\lambda\,{\rm id})).$$
On the other hand, by the assumption, we have $\pm\beta(g_{\ast}^{-1}
\overline v)\notin{\rm Spec}(A_{\overline v}^{\bf c}
\vert_{g_{\ast}\mathfrak q_{\beta}^{\bf c}})$ for each 
$\beta\in\overline{\triangle}_+$.  
Therefore, it follows from Lemma 5.3, Propositions 5.6 and 5.7 of [17] that 
there exists a pseudo-orthonormal base of $(T_u\widetilde M)^{\bf c}$ 
consisting of eigenvectors of $\widetilde A_v^{\bf c}$.  Therefore 
$\widetilde M$ is proper complex isoparametric, that is, $M$ is proper 
complex equifocal.  \hspace{1.5truecm}q.e.d.

\section{Proof of Theorems A, C and E} 
In this section, we shall prove Theorems A, C and E.  
First we prove Theroem A in terms of Propositions 3.2, 4.1 and 5.2.  

\vspace{0.5truecm}

\noindent
{\it Proof of Theorem A.} Since $T_{eH}(H'(eH))=\mathfrak q\cap\mathfrak h'$ 
and $\mathfrak q\cap\mathfrak h'$ is a non-degenerate subsapce of 
$\mathfrak q$, we see that $H'(eH)$ is a pseudo-Riemannian submanifold.  Since 
$\sigma\circ\sigma'=\sigma'\circ\sigma$, we can show that $H'(eH)$ is 
a reflective submanifold by imitating the first-half part 
of the proof of Lemma 4.2 in [19].  Thus the first-half part of the statement 
(i) is shown.  
Furthermore, by imitating the second-half part of the proof of Lemma 4.2 in 
[19], we can show the second-half part of the statement (i).  
In the sequel, we shall 
show the statement (ii).  Let $M$ be a principal orbit of 
the $H'$-action as in the statement (ii).  
For simplicity, set $x:=\exp_G(w)H$ and $g:=\exp_G(w)$, where $w$ is as in 
the statement (ii).  
By imitating the second-half part of the proof of Lemma 4.2 in [17], it is 
shown that $M$ is a partial tube over $H'(eH)$ and $M\cap\Sigma_{eH}$ is 
an orbit of the isotropy action of the symmetric space 
$\Sigma_{eH}(\cong L/H\cap H')$.  
Since $M$ is a principal orbit, $M\cap \Sigma_{eH}$ is 
a principal orbit of the isotropy action.  
Hence, since ${\rm ad}(w)\vert_{\mathfrak l}$ is semi-simple, 
$\mathfrak b:=g_{\ast}^{-1}T^{\perp}_xM$ is a Cartan subspace of 
$\mathfrak q\cap\mathfrak q'$ by (i) of Proposition 3.2.  
Take a Cartan subspace $\mathfrak a$ of $\mathfrak q$ containing 
$\mathfrak b$.  Let $\mathfrak q^{\bf c}=\mathfrak a^{\bf c}+
\sum_{\alpha\in\triangle_+}\mathfrak q_{\alpha}^{\bf c}$ be the root space 
decomposition with respect to $\mathfrak a^{\bf c}$.  Set 
$\triangle_{\mathfrak b^{\bf c}}:=\{\alpha\vert_{\mathfrak b^{\bf c}}\,\vert\,
\alpha\in\triangle\,\,{\rm s.t.}\,\,\alpha\vert_{\mathfrak b^{\bf c}}\not=0\}$ 
and $\mathfrak q_{\beta}^{\bf c}:=\sum_{\alpha\in\triangle\,\,{\rm s.t.}\,\,
\alpha\vert_{\mathfrak b^{\bf c}}=\beta}\mathfrak q_{\alpha}^{\bf c}$ 
($\beta\in\triangle_{\mathfrak b^{\bf c}}$).  Then we have $\mathfrak q^{\bf c}
=\mathfrak z_{\mathfrak q^{\bf c}}(\mathfrak b^{\bf c})+
\sum_{\beta\in(\triangle_{\mathfrak b^{\bf c}})_+}\mathfrak q_{\beta}^{\bf c}$, where $(\triangle_{\mathfrak b^{\bf c}})_+$ is the positive root system 
under some lexicographical ordering.  
Also, since $\mathfrak q^{\bf c}\cap{\mathfrak h'}^{\bf c}$ and 
$\mathfrak q^{\bf c}\cap{\mathfrak q'}^{\bf c}$ are ${\rm ad}(b)^2$-invariant 
for any $b\in\mathfrak b^{\bf c}$, we have 
$\mathfrak q^{\bf c}\cap{{\mathfrak h}'}^{\bf c}=
\mathfrak z_{\mathfrak q^{\bf c}}(\mathfrak b^{\bf c})\cap
{{\mathfrak h}'}^{\bf c}+\sum_{\beta\in(\triangle_{\mathfrak b^{\bf c}})_+}
(\mathfrak q_{\beta}^{\bf c}\cap{{\mathfrak h}'}^{\bf c})$ and 
$\mathfrak q^{\bf c}\cap{{\mathfrak q}'}^{\bf c}=
\mathfrak b^{\bf c}+\sum_{\beta\in(\triangle_{\mathfrak b^{\bf c}})_+}
(\mathfrak q_{\beta}^{\bf c}\cap{{\mathfrak q}'}^{\bf c})$.  
Hence we have 
$$
(T_xM)^{\bf c}=g_{\ast}^{\bf c}
(\mathfrak z_{\mathfrak q^{\bf c}}
(\mathfrak b^{\bf c})\cap{{\mathfrak h}'}^{\bf c})+
\sum_{\beta\in(\triangle_{\mathfrak b^{\bf c}})_+}
\left(
g_{\ast}^{\bf c}(\mathfrak q_{\beta}^{\bf c}\cap{{\mathfrak h}'}^{\bf c})
+g_{\ast}^{\bf c}(\mathfrak q_{\beta}^{\bf c}\cap{{\mathfrak q}'}^{\bf c})
\right),
$$
$$(T_{eH}(H'(eH)))^{\bf c}=
\mathfrak z_{\mathfrak q^{\bf c}}(\mathfrak b^{\bf c})\cap
{{\mathfrak h}'}^{\bf c}+\sum_{\beta\in(\triangle_{\mathfrak b^{\bf c}})_+}
(\mathfrak q_{\beta}^{\bf c}\cap{{\mathfrak h}'}^{\bf c})$$
and 
$$(T_x(M\cap\Sigma_{eH}))^{\bf c}=
\sum_{\beta\in(\triangle_{\mathfrak b^{\bf c}})_+}
g_{\ast}^{\bf c}(\mathfrak q_{\beta}^{\bf c}\cap{{\mathfrak q}'}^{\bf c}).$$
Also we have $T^{\perp}_xM=g_{\ast}\mathfrak b$.  
Take $v\in T^{\perp}_xM=g_{\ast}\mathfrak b$.  
It is clear that $R(\cdot,v)v$ is semi-simple.  
Since $H'(eH)$ is totally geodesic, it follows from (ii) of Proposition 4.1 
and $(4.1)$ that 
$A^{\bf c}_v\widetilde X_w=0\quad(X\in\mathfrak z_{\mathfrak q^{\bf c}}
(\mathfrak b^{\bf c})\cap{\mathfrak h'}^{\bf c})$ 
and 
$$A^{\bf c}_v\widetilde X_w=\sqrt{-1}\beta(g_{\ast}^{-1}v)\tan
(\sqrt{-1}\beta(w))\widetilde X_w
\quad(X\in\mathfrak q_{\beta}^{\bf c}\cap{\mathfrak h'}^{\bf c}\,\,
(\beta\in(\triangle_{\mathfrak b^{\bf c}})_+)).\leqno{(6.1)}$$
Also, since $M\cap\Sigma_{eH}$ is a principal orbit of 
the isotropy action of $\Sigma_{eH}(\cong L/H\cap K)$, it follows from 
Proposition 3.2 and (i) of Proposition 4.1 that 
$$A_v^{\bf c}Y=-\frac{\sqrt{-1}\beta(g_{\ast}^{-1}v)}{\tan(\sqrt{-1}\beta(w))}Y
\quad(Y\in g_{\ast}(\mathfrak q_{\beta}^{\bf c}\cap{\mathfrak q'}^{\bf c}))
\leqno{(6.2)}$$
up to constant-multiple, where we note that the induced metric on 
$\Sigma_{eH}(=L/H\cap K)$ is homothetic to the metric induced from 
the Killing form of $\mathfrak l$.  
Thus $A^{\bf c}_v$ is diagonalizable, that is, it is semi-simple and we have 

\noindent
$[A_v^{\bf c},R^{\bf c}(\cdot,v)v\vert_{(T_xM)^{\bf c}}]=0$ and hence 
$[A_v,R(\cdot,v)v\vert_{T_xM}]=0$.  
Therefore $M$ is curvature-adapted.  
Next we shall show that $M$ is proper complex equifocal.  
Since $g_{\ast}^{-1}T^{\perp}_xM$ is a Cartan subspace of 
$\mathfrak q\cap\mathfrak q'$ for each $x(=gH)\in M$, 
$M$ has flat section.  Since $M$ is a principal orbit of 
the $H'$-action, each normal vector of $M$ extend to an 
$H'$-equivariant normal 
vector field, which is parallel with respect to the normal connection of $M$ 
because $M$ has flat section.  
From this fact, it follows that the normal holonomy group of $M$ is trivial.  
Furthermore, it follows from the homogeneity of $M$ that $M$ is complex 
equifocal.  From $(6.1)$ and $(6.2)$, we have 
${\rm Spec}(A_v^{\bf c}\vert_{g_{\ast}^{\bf c}\mathfrak q_{\beta}^{\bf c}})
\subset\{\sqrt{-1}\beta(g_{\ast}^{-1}v)\tan(\sqrt{-1}\beta(w)),\,\,
-\frac{\sqrt{-1}\beta(g_{\ast}^{-1}v)}
{\tan(\sqrt{-1}\beta(w))}\}$ ($\beta\in(\triangle_{\mathfrak b^{\bf c}})_+$), 
that is, 

\noindent
$\pm\beta(g_{\ast}^{-1}v)\notin{\rm Spec}\,A^{\bf c}_v
\vert_{g_{\ast}^{\bf c}\mathfrak q^{\bf c}_{\beta}}$.  
Therefore, it follows from Proposition 5.2 that 
$M$ is proper complex equifocal.  
Furthermore it follows from the result of [23] stated in Introduction that 
$M$ is an isoparametric submanifold with flat section.  
This completes the proof.  \hspace{1.5truecm}q.e.d.

\vspace{0.5truecm}

Next we prove Theorem C.  

\vspace{0.5truecm}

\noindent
{\it Proof of Theorem C.} According to Theorem A, we have only to show that 
$K(eH)$ has no focal point and that, for any normal vector $v$ of $M_i$, 
$R(\cdot,v)v\vert_{T_xM_i}$ and $A_v$ are 
diagonalizable.  Let $\mathfrak g=\mathfrak f+\mathfrak p$ be 
the Cartan decomposition of 
$\mathfrak g$ associated with $\theta$.  
Take an arbitrary normal vector $v$ of $K(eH)$ at $eH$.  Take a maximal 
abelian subspace $\mathfrak b$ of $\mathfrak q\cap\mathfrak p$ containing $v$ 
and a Cartan subspace $\mathfrak a$ of $\mathfrak q$ containing 
$\mathfrak b$.  Let $\mathfrak q^{\bf c}=\mathfrak a^{\bf c}+
\sum_{\alpha\in\triangle_+}\mathfrak q_{\alpha}^{\bf c}$ be the root space 
decomposition of $\mathfrak q^{\bf c}$ with respect to $\mathfrak a^{\bf c}$.  
Let $\triangle_{\mathfrak b}:=
\{\alpha\vert_{\mathfrak b}\,\vert\,\alpha\in\triangle\,\,{\rm s.t.}\,\,
\alpha\vert_{\mathfrak b}\not=0\}$ and $\mathfrak q_{\beta}
:=(\sum_{\alpha\in\triangle\,\,{\rm s.t.}\,\,\alpha\vert_{\mathfrak b}=\beta}
\mathfrak q_{\alpha}^{\bf c})\cap\mathfrak q$ 
($\beta\in\triangle_{\mathfrak b}$).  
Since $\mathfrak b\subset\mathfrak p$, we have $\beta(\mathfrak b)\subset
{\bf R}$ ($\beta\in\triangle_{\mathfrak b}$) (see Lemma 3.1) and hence 
$\mathfrak q=\mathfrak z_{\mathfrak q}(\mathfrak b)+
\sum_{\beta\in(\triangle_{\mathfrak b})_+}\mathfrak q_{\beta}$.  Furthermore, 
since ${\rm ad}(b)^2(\mathfrak q\cap\mathfrak f)\subset\mathfrak q\cap
\mathfrak f$ for any $b\in\mathfrak b$, we have $\mathfrak q\cap\mathfrak f=
\mathfrak z_{\mathfrak q}(\mathfrak b)\cap\mathfrak f+
\sum_{\beta\in(\triangle_{\mathfrak b})_+}(\mathfrak q_{\beta}\cap
\mathfrak f)$.  Let $X\in\mathfrak q_{\beta}\cap\mathfrak f$ ($\beta\in
(\triangle_{\mathfrak b})_+$), $Y$ be the strongly $K(eH)$-Jacobi field along 
$\gamma_v$ with $Y(0)=X$.  Since $K(eH)$ is totally geodesic, we have 
$Y(s)=\cosh(s\beta(v))P_{\gamma_v\vert_{[0,s]}}(X)$.  Since $\beta(v)$ is a 
real number, $Y$ has no zero point.  Also any strongly $K(eH)$-Jacobi field 
$\widehat Y$ along $\gamma_v$ with $\widehat Y(0)\in\mathfrak z_{\mathfrak q}
(\mathfrak b)\cap\mathfrak f$ is expressed as $\widehat Y(s)=
P_{\gamma_v\vert_{[0,s]}}(\widehat Y(0))$ and hence it has no zero point.  
On the other hand, since $K(eH)$ is reflective and hence it has section, 
any non-strongly $K(eH)$-Jacobi field along $\gamma_v$ has no zero point.  
After all there exists no focal point of $K(eH)$ along $\gamma_v$.  
From the arbitrariness of $v$, it follows 
that $K(eH)$ has no focal point.  For convenience, set $H_1:=K,\,H_2:=L,\,
\mathfrak h_1:=\mathfrak f,\,\mathfrak h_2:=\mathfrak l,\,
\mathfrak q_1:=\mathfrak p$ and $\mathfrak q_2:=\mathfrak f\cap\mathfrak q
+\mathfrak p\cap\mathfrak h$.  Let $M_1$ (resp. $M_2$) be a 
principal orbit of the $H_1$-action (resp. the $H_2$-action) 
through $x_1=\exp_G(w_1)H\in H_2(eH)$ ($w_1\in \mathfrak q\cap\mathfrak q_1$) 
(resp. $x_2=\exp_G(w_2)H\in H_1(eH)\setminus F$ ($w_2\in \mathfrak q\cap
\mathfrak q_2$)).  
Set $g_i:=\exp_G(w_i)$ ($i=1,2$).  Since 
$\mathfrak b_1:=g_{1\ast}^{-1}(T^{\perp}_{x_1}M_1)$ and 
$\mathfrak b_2:=g_{2\ast}^{-1}(T^{\perp}_{x_2}M_2)$ are maximal abelian 
subspaces of $\mathfrak q\cap\mathfrak p$ and $\mathfrak q\cap\mathfrak f$, 
respectively, they are maximal split abelian subspaces of $\mathfrak q$.  
Hence we have 
the root space decomposition $\mathfrak q=\mathfrak z_{\mathfrak q}
(\mathfrak b_i)+\sum_{\beta\in\triangle^i_+}\mathfrak q_{\beta}$ of 
$\mathfrak q$ with respect to $\mathfrak b_i$ ($i=1,2$), where 
$\mathfrak q_{\beta}:=\{X\in\mathfrak q\,\vert\,{\rm ad}(b)^2(X)=(-1)^{i-1}
\beta(b)^2X\,\,(\forall\,\,b\in\mathfrak b_i)\}$ by Lemma 3.1 
and $\triangle^i_+$ is 
the positive root system of $\triangle^i:=\{\beta\in\mathfrak b_i^{\ast}\,
\vert\,\mathfrak q_{\beta}\not=\{0\}\}$ with respect to a lexicographical 
ordering.  Also, it is easy to show that $\mathfrak q\cap\mathfrak h_i
=\mathfrak z_{\mathfrak q}(\mathfrak b_i)\cap\mathfrak h_i+
\sum_{\beta\in\triangle^i_+}(\mathfrak q_{\beta}\cap\mathfrak h_i)$ and 
$\mathfrak q\cap\mathfrak q_i=\mathfrak b_i+
\sum_{\beta\in\triangle^i_+}(\mathfrak q_{\beta}\cap\mathfrak q_i)$, 
where $i=1,2$.  Hence we have 
$$
T_{x_i}M_i=g_{i\ast}(\mathfrak z_{\mathfrak q}
(\mathfrak b_i)\cap\mathfrak h_i)+
\sum_{\beta\in\triangle^i_+}
\left(g_{i\ast}(\mathfrak q_{\beta}\cap\mathfrak h_i)
+g_{i\ast}(\mathfrak q_{\beta}\cap\mathfrak q_i)\right),
$$
$$T_{eH}(H_i(eH))=
\mathfrak z_{\mathfrak q}(\mathfrak b_i)\cap\mathfrak h_i
+\sum_{\beta\in\triangle^i_+}(\mathfrak q_{\beta}\cap\mathfrak h_i)$$
and 
$$T_{x_i}(M_i\cap\Sigma^i_{eH})=
\sum_{\beta\in\triangle^i_+}
g_{i\ast}(\mathfrak q_{\beta}\cap\mathfrak q_i),$$
where $\Sigma^i_{eH}$ is the section of $H_i(eH)$ through $eH$.  
Take $v_i\in T^{\perp}_{x_i}M_i=g_{i\ast}\mathfrak b_i$.  
It is clear that $R(\cdot,v_i)v_i$ is diagonalizable.  
Denote by $A^i$ the shape tensor of $M_i$.  
By using Propositions 3.2, 4.1 and $(4.1)$, we can show 
$A^i_{v_i}\widetilde X_{w_i}=0\quad(X\in\mathfrak z_{\mathfrak q}
(\mathfrak b_i)\cap\mathfrak h_i)$, 
$$A^i_{v_i}\widetilde X_{w_i}=\sqrt{-1}^i\beta(g_{i\ast}^{-1}v_i)
\tan(\sqrt{-1}^i\beta(w_i))\widetilde X_{w_i}
\quad(X\in\mathfrak q_{\beta}\cap\mathfrak h_i\,\,
(\beta\in\triangle^i_+))$$
and 
$$A^i_{v_i}Y=-\frac{\sqrt{-1}^i\beta(g_{i\ast}^{-1}v_i)}
{\tan(\sqrt{-1}^i\beta(w_i))}Y\quad(Y\in 
g_{i\ast}(\mathfrak q_{\beta}\cap\mathfrak q_i)\,(\beta\in\triangle^i_+)).$$
Thus $A^i_{v_i}$ is diagonalizable.  
This completes the proof.  \hspace{1.5truecm}q.e.d.

\vspace{0.5truecm}

Next we shall prove Theorem E.  
By imitating the proof of Lemma 2.1 of [21], we can show the following fact.  

\vspace{0.5truecm}

\noindent
{\bf Lemma 6.1.} {\sl Let $G(=(G\times G)/\triangle G)$ be a semi-simple Lie 
group equipped with the bi-invariant pseudo-Riemannian metric induced from 
the Killing form of $\mathfrak g+\mathfrak g$, $H'$ be a closed subgroup of 
$G\times G$ and $\mathfrak a$ be an abelian subspace of the normal space 
$T^{\perp}_e(H'\cdot e)$ of $H'\cdot e$.  Set $\Sigma:=\exp_G(\mathfrak a)$.  
Then all $H'$-orbits through $\Sigma$ meet $\Sigma$ orthogonally.}

\vspace{0.5truecm}

By using this lemma and imitating the proof of Lemma 2.2 of [21], 
we can show  tyhe following fact.  

\vspace{0.5truecm}

\noindent
{\bf Lemma 6.2.} {\sl Let $G/H$ be a semi-simple pseudo-Riemannian symmetric 
space, $H'$ be a closed subgroup of $G$ and $\mathfrak a$ be an abelian 
subspace of the normal space $T^{\perp}_{eH}H(eH)$ of $H'(eH)$.  
Set $\Sigma:={\rm Exp}(\mathfrak a)$.  Then all $H'$-orbits through $\Sigma$ 
meet $\Sigma$ orthogonally.}

\vspace{0.5truecm}

By using this lemma, we prove Theorem E.  

\vspace{0.5truecm}

\noindent
{\it Proof of Theorem E.} Let $M,\,F$ and $G/H$ be as in the statement of 
Theorem E.  Without loss of generality, we may assume that $G$ is simply 
connected.  Since $M$ is homogeneous, there exists a closed subgroup $H_1$ of 
$G$ having $M$ as an orbit.  Without loss of generality, we may assume that 
$H_1(eH)=M$.  Set $\Sigma:={\rm Exp}(T^{\perp}_{eH}M)$.  Since $M$ has flat 
section, that is, $T^{\perp}_{eH}M$ is abelian, it follows from Lemma 6.2 that 
all $H_1$-orbits through $\Sigma$ meet $\Sigma$ orthogonally.  Hence their 
dimensions are lower than ${\rm dim}\,M+1$.  This fact together with 
${\rm dim}\,M+{\rm dim}\,\Sigma={\rm dim}\,G/H$ implies that all $H_1$-orbits 
through $W$ are of the same dimension as ${\rm dim}\,M$ for some neighborhood 
$W$ of $eH$ in $\Sigma$ and they are principal orbits.  Set $U:=H_1\cdot W$, 
which is an open set of $G/H$.  Fix $g_0H\in F$.  Set $H_2:=g_0^{-1}H_1g_0, 
\mathfrak t:=T_{eH}g_0^{-1}F$ and $\mathfrak t^{\perp}:=T_{eH}^{\perp}g_0^{-1}
F$.  Furthermore set $\mathfrak h':=\mathfrak n_{\mathfrak h}(\mathfrak t)
+\mathfrak t$ and $\mathfrak q':=(\mathfrak h\ominus\mathfrak n_{\mathfrak h}
(\mathfrak t))+\mathfrak t^{\perp}$.  Since $\mathfrak n_{\mathfrak h}
(\mathfrak t)$ is a non-degenerate subspace of $\mathfrak h$ by the 
assumption, we have $\mathfrak g=\mathfrak h'\oplus\mathfrak q'$ 
(orthogonal direct sum).  Since $F$ is a reflective by the assumption, 
$\mathfrak t$ and $\mathfrak t^{\perp}$ are Lie triple systems.  
By using this fact, we can show $[\mathfrak h',\mathfrak h']
\subset\mathfrak h',\,[\mathfrak h',\mathfrak q']\subset\mathfrak q'$ and 
$[\mathfrak q',\mathfrak q']\subset\mathfrak h'$.  Thus the connected subgroup 
$H'$ of $G$ having $\mathfrak h'$ as its Lie algebra is symmetric, where we 
use the simply connectedness of $G$.  That is, the $H'$-action on $G/H$ is 
a Hermann type action.  Easily we can show $T_e((H_2\times H)\cdot e)
={\rm pr}_{\mathfrak q}(\mathfrak h_2)+\mathfrak h$ and 
$T_e((H'\times H)\cdot e)={\rm pr}_{\mathfrak q}(\mathfrak h')+\mathfrak h
=\mathfrak t+\mathfrak h$, where ${\rm pr}_{\mathfrak q}$ is the orthogonal 
projection of $\mathfrak g$ onto $\mathfrak q$ and $\mathfrak h_2:=
{\rm Lie}\,H_2$.  Since $\pi^{-1}(H_2(eH))=(H_2\times H)\cdot e$, we have 
$T_e(H_2(eH))={\rm pr}_{\mathfrak q}(T_e((H_2\times H)\cdot e))
={\rm pr}_{\mathfrak q}(\mathfrak h_2)$, that is, ${\rm pr}_{\mathfrak q}
(\mathfrak h_2)=\mathfrak t$.  Hence we have $T_e((H'\times H)\cdot e)
=T_e((H_2\times H)\cdot e))$, which implies $(H'\times H)\cdot e=(H_2\times H)
\cdot e$.  Therefore we have $H'(eH)=H_2(eH)$.  Set $\Sigma':={\rm Exp}
(T^{\perp}_{g_0^{-1}H}(g_0^{-1}M))$, which passes through $eH$.  Set 
$\mathfrak a':=T_{eH}\Sigma'$, which is abelian.  Since 
$T^{\perp}_{eH}(H'(eH))=T^{\perp}_{eH}(H_2(eH))$ includes $\mathfrak a'$, it 
follows from Lemma 6.2 that all $H'$-orbits and all $H_2$-orbits through 
$\Sigma'$ meet $\Sigma'$ orthogonally.  Since all $H_2$-orbits through 
$g_0^{-1}W(\subset \Sigma')$ are principal and hence 
$T^{\perp}_{gH}(H_2(gH))=T_{gH}\Sigma'$ 
for all $gH\in g_0^{-1}W$, we have $T_{gH}(H'(gH))\subset T_{gH}(H_2(gH))$ 
for all $gH\in g_0^{-1}W$.  On the other hand, we have 
$[{\rm pr}_{\mathfrak h}(\mathfrak h_2),\mathfrak t]={\rm pr}_{\mathfrak q}
([\mathfrak h_2,\mathfrak t])\subset{\rm pr}_{\mathfrak q}
(T_e((H_2\times H)\cdot e))=T_{eH}((H_2(eH))=\mathfrak t$, that is, 
${\rm pr}_{\mathfrak h}(\mathfrak h_2)\subset\mathfrak n_{\mathfrak h}
(\mathfrak t)$, where ${\rm pr}_{\mathfrak h}$ is the orthogonal projection of 
$\mathfrak g$ onto $\mathfrak h$.  Hence we have 
$\mathfrak h_2\subset{\rm pr}_{\mathfrak h}(\mathfrak h_2)
+{\rm pr}_{\mathfrak q}(\mathfrak h_2)\subset\mathfrak h'$, that is, 
$H_2\subset H'$.  Therefore we see that $H'(gH)=H_2(gH)$ for all 
$gH\in g_0^{-1}W$.  In particular, $g_0^{-1}M$ is a principal orbit of the 
$H'$-action.  Hence $M$ is a principal orbit of the Hermann type action 
$g_0H'g_0^{-1}$.  
This completes the proof.  \hspace{1.5truecm}q.e.d.

\section{Cohomogeneities of special Hermann type actions} 
In this section, we shall list up the cohomogeneities of the $K$-action and 
the $L$-action as in Theorem C on irreducible (semi-simple) pseudo-Riemannian 
symmetric spaces $G/H$ in terms of the fact that the cohomogeneity of 
the $K$-action (resp. $L$-action) is equal to the rank of $L/H\cap K$ 
(resp. $K/H\cap K$).  
In Tables $1\sim5$, $A\cdot B$ denotes $A\times B/\Pi$, where 
$\Pi$ is the discrete center of $A\times B$.  
The symbol $\widetilde{SO_0(1,8)}$ in Table 6 denotes the universal covering 
of $SO_0(1,8)$ and the symbol $\alpha$ in Table 6 denotes an outer 
automorphism of $G_2^2$.  

\newpage


$$
\begin{tabular}{|c|c|c|}
\hline
{\scriptsize {\scriptsize$G/H$}} & 
{\scriptsize$K$} & {\scriptsize$L$}\\
\cline{2-2}\cline{3-3}
 & {\scriptsize ${\rm cohom}_K$} & {\scriptsize ${\rm cohom}_L$}\\
\hline
{\scriptsize$SL(n,{\bf R})/SO_0(p,n-p)\,\,(p\leq\frac n2)$} & 
{\scriptsize$SO(n)$} & {\scriptsize$(SL(p,{\bf R})\times SL(n-p,{\bf R}))$}\\
&&{\scriptsize$\cdot{\bf R}_{\ast}$}\\
\cline{2-2}\cline{3-3}
&{\scriptsize$n-1$} & {\scriptsize$p$}\\
\hline
{\scriptsize$SL(n,{\bf R})/(SL(p,{\bf R})\times SL(n-p,{\bf R}))
\cdot{\bf R}_{\ast}$} & {\scriptsize$SO(n)$} 
& {\scriptsize$SO_0(p,n-p)$}\\
\cline{2-2}\cline{3-3}
{\scriptsize$(p\leq\frac n2)$}&{\scriptsize$p$} & {\scriptsize$p$}\\
\hline
{\scriptsize$SL(2n,{\bf R})/Sp(n,{\bf R})$} & {\scriptsize$SO(2n)$} & 
{\scriptsize$SL(n,{\bf C})\cdot U(1)$}\\
\hline
\end{tabular}$$

\centerline{{\bf Table 1.}}

\newpage

$$\begin{tabular}{|c|c|c|}
\hline
{\scriptsize {\scriptsize$G/H$}} & 
{\scriptsize$K$} & {\scriptsize$L$}\\
\cline{2-2}\cline{3-3}
 & {\scriptsize ${\rm cohom}_K$} & {\scriptsize ${\rm cohom}_L$}\\
\hline
{\scriptsize$SL(n,{\bf C})/SO(n,{\bf C})$} & {\scriptsize$SU(n)$} & 
{\scriptsize$SL(n,{\bf R})$}\\
\cline{2-2}\cline{3-3}
 & {\scriptsize$n-1$} & {\scriptsize$n-1$}\\
\hline
{\scriptsize$SL(n,{\bf C})/SL(n,{\bf R})$} & {\scriptsize$SU(n)$} & 
{\scriptsize$SO(n,{\bf C})$}\\
\cline{2-2}\cline{3-3}
 & {\scriptsize$[\frac n2]$} & {\scriptsize$n-1$}\\
\hline
{\scriptsize$SL(n,{\bf C})/(SL(p,{\bf C})\times SL(n-p,{\bf C})\times U(1))$} 
& {\scriptsize$SU(n)$} & {\scriptsize$SU(p,n-p)$}\\
\cline{2-2}\cline{3-3}
{\scriptsize $(p\leq \frac n2)$} & {\scriptsize $p$} & {\scriptsize $p$}\\
\hline
\end{tabular}$$

\centerline{{\bf Table 2.}}

\newpage

$$\begin{tabular}{|c|c|c|}
\hline
{\scriptsize {\scriptsize$G/H$}} & 
{\scriptsize$K$} & {\scriptsize$L$}\\
\cline{2-2}\cline{3-3}
 & {\scriptsize ${\rm cohom}_K$} & {\scriptsize ${\rm cohom}_L$}\\
\hline
{\scriptsize$Sp(n,{\bf R})/SU(p,n-p)\cdot U(1)\,\,(p\leq\frac n2)$} & 
{\scriptsize$U(n)$} & {\scriptsize$Sp(p,{\bf R})$}\\
&&{\scriptsize$\times Sp(n-p,{\bf R})$}\\
\cline{2-2}\cline{3-3}
& {\scriptsize$n$} & {\scriptsize$p$}\\
\hline
{\scriptsize$Sp(n,{\bf R})/Sp(p,{\bf R})\times Sp(n-p,{\bf R})$} & 
{\scriptsize$U(n)$} & {\scriptsize$SU(p,n-p)\cdot U(1)$}\\
\cline{2-2}\cline{3-3}
{\scriptsize$(p\leq\frac n2)$} & {\scriptsize$p$} & {\scriptsize$p$}\\
\hline
{\scriptsize$Sp(n,{\bf R})/SL(n,{\bf R})\cdot U(1)$} & {\scriptsize$U(n)$} & 
{\scriptsize$SL(n,{\bf R})\cdot U(1)$}\\
\hline
\end{tabular}$$

\centerline{{\bf Table 3.}}

\newpage

$$\begin{tabular}{|c|c|c|c|c|}
\hline
{\scriptsize$G/H$} & {\scriptsize$K$} & {\scriptsize$L$} 
& {\scriptsize${\rm cohom}_K$} & {\scriptsize${\rm cohom}_L$}\\
\hline
{\scriptsize$E_6^6/Sp(4,{\bf R})$} & {\scriptsize$Sp(4)/\{\pm1\}$} 
& {\scriptsize$SL(6,{\bf R})\times SL(2,{\bf R})$} 
& {\scriptsize$6$} & {\scriptsize$4$}\\
\hline
{\scriptsize$E_6^6/SL(6,{\bf R})\times SL(2,{\bf R})$} & 
{\scriptsize$Sp(4)/\{\pm1\}$} & {\scriptsize$Sp(4,{\bf R})$} 
& {\scriptsize$4$} & {\scriptsize$4$}\\
\hline
{\scriptsize$E_6^6/Sp(2,2)$} & {\scriptsize$Sp(4)/\{\pm1\}$} & 
{\scriptsize$SO_0(5,5)\cdot{\bf R}$} & {\scriptsize$6$} & {\scriptsize$2$}\\
\hline
{\scriptsize$E_6^6/SO_0(5,5)\cdot{\bf R}$} & {\scriptsize$Sp(4)/\{\pm1\}$} & 
{\scriptsize$Sp(2,2)$} & {\scriptsize$2$} & {\scriptsize$2$}\\
\hline
{\scriptsize$E_6^6/SU^{\ast}(6)\cdot SU(2)$} & {\scriptsize$Sp(4)/\{\pm1\}$} & 
{\scriptsize$F^4_4$} & {\scriptsize$4$} & {\scriptsize$1$}\\
\hline
{\scriptsize$E_6^6/F^4_4$} & {\scriptsize$Sp(4)/\{\pm1\}$} & 
{\scriptsize$SU^{\ast}(6)\cdot SU(2)$} & {\scriptsize$2$} & {\scriptsize$1$}\\
\hline
{\scriptsize$E_6^2/Sp(1,3)$} & {\scriptsize$SU(6)\cdot SU(2)$} 
& {\scriptsize$F_4^4$} & {\scriptsize$4$} & {\scriptsize$2$}\\
\hline
\end{tabular}$$

\centerline{{\bf Table 4.}}

\newpage


$$\begin{tabular}{|c|c|c|c|c|}
\hline
{\scriptsize$G/H$} & {\scriptsize$K$} & {\scriptsize$L$} 
& {\scriptsize${\rm cohom}_K$} & {\scriptsize${\rm cohom}_L$}\\
\hline
{\scriptsize$E_7^7/SL(8,{\bf R})$} & {\scriptsize$SU(8)/\{\pm1\}$} & 
{\scriptsize$SL(8,{\bf R})$} & {\scriptsize$7$} & {\scriptsize$7$}\\
\hline
{\scriptsize$E_7^7/SU^{\ast}(8)$} & {\scriptsize$SU(8)/\{\pm1\}$} & 
{\scriptsize$E_6^6\cdot U(1)$} & {\scriptsize$7$} & {\scriptsize$3$}\\
\hline
{\scriptsize$E_7^7/E_6^6\cdot U(1)$} & {\scriptsize$SU(8)/\{\pm1\}$} & 
{\scriptsize$SU^{\ast}(8)$} & {\scriptsize$3$} & {\scriptsize$3$}\\
\hline
{\scriptsize$E_7^7/SU(4,4)$} & {\scriptsize$SU(8)/\{\pm1\}$} & 
{\scriptsize$SO_0(6,6)\cdot SL(2,{\bf R})$} & {\scriptsize$7$}
 & {\scriptsize$4$}\\
\hline
{\scriptsize$E_7^7/SO_0(6,6)\cdot SL(2,{\bf R})$}
 & {\scriptsize$SU(8)/\{\pm1\}$} & {\scriptsize$SU(4,4)$} & {\scriptsize$4$}
 & {\scriptsize$4$}\\
\hline
{\scriptsize$E_7^7/SO^{\ast}(12)\cdot SU(2)$}
 & {\scriptsize$SU(8)/\{\pm1\}$} & 
{\scriptsize$E_6^2\cdot U(1)$} & {\scriptsize$4$} & {\scriptsize$2$}\\
\hline
{\scriptsize$E_7^7/E_6^2\cdot U(1)$}
 & {\scriptsize$SU(8)/\{\pm1\}$} & 
{\scriptsize$SO^{\ast}(12)\cdot SU(2)$} & {\scriptsize$3$} & {\scriptsize$2$}\\
\hline
\end{tabular}$$

\centerline{{\bf Table 5.}}

\newpage

$$\begin{tabular}{|c|c|c|c|c|}
\hline
{\scriptsize$G/H$} & {\scriptsize$K$} & {\scriptsize$L$} 
& {\scriptsize${\rm cohom}_K$} & {\scriptsize${\rm cohom}_L$}\\
\hline
{\scriptsize$F_4^4/Sp(1,2)\cdot Sp(1)$} & {\scriptsize$Sp(3)\cdot Sp(1)$} 
& {\scriptsize$SO_0(4,5)$} & {\scriptsize$4$} & {\scriptsize$1$}\\
\hline
{\scriptsize$F_4^4/SO_0(4,5)$} & {\scriptsize$Sp(3)\cdot Sp(1)$} 
& {\scriptsize$Sp(1,2)\cdot Sp(1)$} & {\scriptsize$1$} & {\scriptsize$1$}\\
\hline
{\scriptsize$F_4^4/Sp(3,{\bf R})\cdot SL(2,{\bf R})$} 
& {\scriptsize$Sp(3)\cdot Sp(1)$} & 
{\scriptsize$Sp(3,{\bf R})\cdot SL(2,{\bf R})$} & {\scriptsize$4$}
 & {\scriptsize$4$}\\
\hline
{\scriptsize$F^{-20}_4/Sp(1,2)\cdot Sp(1)$} & {\scriptsize$Spin(9)$}
 & {\scriptsize$\widetilde{SO_0(1,8)}$} & {\scriptsize$1$} & {\scriptsize$1$}\\
\hline
{\scriptsize$F^{-20}_4/\widetilde{SO_0(1,8)}$} & {\scriptsize$Spin(9)$}
 & {\scriptsize$Sp(1,2)\cdot Sp(1)$} & {\scriptsize$1$} & {\scriptsize$1$}\\
\hline
{\scriptsize$F_4^{\bf C}/F_4^4$} & {\scriptsize$F_4$} & 
{\scriptsize$Sp(3,{\bf C})\cdot SL(2,{\bf C})$} & {\scriptsize$4$}
 & {\scriptsize$4$}\\
\hline
{\scriptsize$F_4^{\bf C}/Sp(3,{\bf C})\cdot SL(2,{\bf C})$}
 & {\scriptsize$F_4$} & {\scriptsize$F_4^4$} & {\scriptsize$4$}
 & {\scriptsize$4$}\\
\hline
\end{tabular}$$

\centerline{{\bf Table 6.}}

\vspace{0.5truecm}

\centerline{{\bf References}}

\vspace{0.3truecm}

{\small
\noindent
[1] M. Berger, Les espaces sym$\acute e$triques non compacts, 
Ann. Sci. $\acute E$c. Norm. Sup$\acute e$r. III. S$\acute e$r. 

{\bf 74} (1959) 85-177.

\noindent
[2] J. Berndt, Real hypersurfaces with constant principal curvatures 
in complex hyperbolic 

space, J. Reine Angew. Math. {\bf 395} (1989) 132-141. 

\noindent
[3] J. Berndt, Real hypersurfaces in quaternionic space forms, 
J. Reine Angew. Math. {\bf 419} 

(1991) 9-26. 

\noindent
[4] J. Berndt and H. Tamaru, Cohomogeneity one actions on noncompact 
symmetric spaces 

with a totally geodesic singular orbit, 
Tohoku Math. J. {\bf 56} (2004) 163-177.

\noindent
[5] J. Berndt and L. Vanhecke, 
Curvature adapted submanifolds, 
Nihonkai Math. J. {\bf 3} (1992) 

177-185.

\noindent
[6] A. Borowiec, M. Francaviglia and I. Volovich, 
Anti-K$\ddot a$hlerian manifolds, 
Differential geom. 

and Its Appl. {\bf 12} (2000) 281-289.

\noindent
[7] L. Geatti, 
Invariant domains in the complexfication of a noncompact Riemannian 
symmetric 

space, J. of Algebra {\bf 251} (2002) 619-685.

\noindent
[8] L. Geatti, 
Complex extensions of semisimple symmetric spaces, manuscripta math. {\bf 120} 

(2006) 1-25.

\noindent
[9] L. Geatti and C. Gorodski, 
Polar orthogonal representations of real reductive algebraic groups, 

arXiv:math.RT/0801.0574v1.

\noindent
[10] O. Goertsches and G. Thorbergsson, 
On the Geometry of the orbits of Hermann actions, 

Geom. Dedicata {\bf 129} (2007) 101-118.

\noindent
[11] J. Hahn, Isoparametric hypersurfaces in the pseudo-Riemannian 
space forms, Math. Z. {\bf 187} 

(1984) 195-208.  

\noindent
[12] J. Hahn, Isotropy representations of semisimple symmetric spaces 
and homogeneous hyper-

surfaces, J. Math. Soc. Japan {\bf 40} (1988) 271-288.  

\noindent
[13] E. Heintze, X. Liu and C. Olmos, 
Isoparametric submanifolds and a 
Chevalley type rest-

riction theorem, Integrable systems, geometry, and topology, 151-190, 
AMS/IP Stud. Adv. 

Math. 36, Amer. Math. Soc., Providence, RI, 2006.

\noindent
[14] E. Heintze, R.S. Palais, C.L. Terng and G. Thorbergsson, 
Hyperpolar actions on symme-

tric spaces, Geometry, topology and physics for Raoul Bott (ed. S. T. Yau), Conf. Proc. 

Lecture Notes Geom. Topology {\bf 4}, 
Internat. Press, Cambridge, MA, 1995 pp214-245.

\noindent
[15] S. Helgason, 
Differential geometry, Lie groups and symmetric spaces, 
Academic Press, New 

York, 1978.

\noindent
[16] I. Kath, 
Indefinite extrinsic symmetric spaces I, arXiv:math.DG/0809.4713v2.

\noindent
[17] N. Koike, 
Submanifold geometries in a symmetric space of non-compact 
type and a pseudo-

Hilbert space, Kyushu J. Math. {\bf 58} (2004) 167-202.

\noindent
[18] N. Koike, 
Complex equifocal submanifolds and infinite dimensional anti-
Kaehlerian isopara-

metric submanifolds, Tokyo J. Math. {\bf 28} (2005) 201-247.

\noindent
[19] N. Koike, 
Actions of Hermann type and proper complex equifocal submanifolds, 
Osaka J. 

Math. {\bf 42} (2005) 599-611.

\noindent
[20] N. Koike, 
A splitting theorem for proper complex equifocal submanifolds, Tohoku Math. J. 

{\bf 58} (2006) 393-417.

\noindent
[21] N. Koike, Complex hyperpolar actions with a totally geodesic orbit, 
Osaka J. Math. {\bf 44} (2007) 

491-503.

\noindent
[22] N. Koike, A Chevalley type restriction theorem for a proper complex 
equifocal submani-

fold, Kodai Math. J. {\bf 30} (2007) 280-296.

\noindent
[23] N. Koike, The complexifications of pseudo-Riemannian manifolds and 
anti-Kaehler geo-

metry, submitted for publication.

\noindent
[24] A. Kollross, A Classification of hyperpolar and cohomogeneity one 
actions, Trans. Amer. 

Math. Soc. {\bf 354} (2001) 571-612.

\noindent
[25] T. Oshima and T. Matsuki, Orbits on affine symmetric spaces under the 
action of the 

isotropy subgroups, J. Math. Soc. Japan {\bf 32} (1980) 399--414.

\noindent
[26] T. Oshima and J. Sekiguchi, The restricted root system of a semisimple 
symmetric pair, 

Advanced Studies in Pure Math. {\bf 4} (1984), 433--497. 

\noindent
[27] B. O'Neill, Semi-Riemannian Geometry, with Applications to Relativity, 
Academic Press, 

New York, 1983.

\noindent
[28] R.S. Palais and C.L. Terng, Critical point theory and submanifold 
geometry, Lecture Notes 

in Math. {\bf 1353}, Springer, Berlin, 1988.

\noindent
[29] W. Rossmann, 
The structure of semisimple symmetric spaces, Can. J. Math. {\bf 1} 
(1979) 157--180.

\noindent
[30] R. Sz$\ddot{{{\rm o}}}$ke, Complex structures on tangent 
bundles of Riemannian manifolds, Math. Ann. {\bf 291} 

(1991) 409-428.

\noindent
[31] R. Sz$\ddot{{{\rm o}}}$ke, Automorphisms of certain Stein 
manifolds, Math. Z. {\bf 219} (1995) 357-385.

\noindent
[32] R. Sz$\ddot{{{\rm o}}}$ke, Adapted complex structures and 
geometric quantization, Nagoya Math. J. {\bf 154} 

(1999) 171-183.

\noindent
[33] R. Sz$\ddot{{{\rm o}}}$ke, Involutive structures on the 
tangent bundle of symmetric spaces, 
Math. Ann. {\bf 319} 

(2001), 319--348.

\noindent
[34] R. Sz$\ddot{{{\rm o}}}$ke, Canonical complex structures associated to 
connections and complexifications of 

Lie groups, Math. Ann. {\bf 329} (2004), 553--591.

\noindent
[35] C.L. Terng, 
Isoparametric submanifolds and their Coxeter groups, 
J. Differential Geometry 

{\bf 21} (1985) 79-107.

\noindent
[36] C.L. Terng and G. Thorbergsson, 
Submanifold geometry in symmetric spaces, J. Differential 

Geometry {\bf 42} (1995) 665-718.

\noindent
[37] D. T$\ddot o$ben, Singular Riemannian foliations on 
nonpositively curved manifolds, Math. Z. {\bf 255} 

(2007) 427-436.

\noindent
[38] L. Verh$\acute o$czki, Shape operators 
of orbits of isotropy subgroups in Riemannian symmetric spa-

ces of the compact type, Beitr$\ddot a$ge zur 
Algebra und Geometrie {\bf 36} (1995) 155-170.
}

\vspace{0.5truecm}

\rightline{Department of Mathematics, Faculty of Science, }
\rightline{Tokyo University of Science}
\rightline{26 Wakamiya Shinjuku-ku,}
\rightline{Tokyo 162-8601, Japan}
\rightline{(e-mail: koike@ma.kagu.tus.ac.jp)}

\end{document}